\documentclass{amsart}

\usepackage{booktabs}
\usepackage{amsmath}
\usepackage{amssymb}
\usepackage{float}
\usepackage{graphicx}
\usepackage{caption}
\usepackage[ruled,linesnumbered]{algorithm2e}
\usepackage{latexsym,enumerate}
\usepackage{mathtools}
\usepackage{fancyhdr}
\usepackage{xspace}
\usepackage{graphicx}
\usepackage{subfigure}
\usepackage{underscore}
\usepackage{hyperref}	
\usepackage{tikz}
\usetikzlibrary{positioning, shapes.geometric}
\hypersetup{colorlinks,
	linkcolor=blue,%
	citecolor=blue}

\usepackage{mathptmx}
\usepackage{mathtools}
\usepackage{graphicx,epstopdf}
\usepackage{bm}
\usepackage{latexsym,enumerate}
\usepackage{mathrsfs}
\usepackage{amsfonts}
\usepackage{amssymb,amsmath}
\usepackage{epsfig}       
\usepackage{subfigure}    
\usepackage{graphicx}
\usepackage{float}
\usepackage{multirow}
\usepackage{array}
\usepackage{color}
\usepackage{listings}
\usepackage{booktabs}
\usepackage{lineno}
\usepackage{fancyhdr}
\newtheorem{theorem}{Theorem}[section]
\newtheorem{lemma}[theorem]{Lemma}

\theoremstyle{definition}
\newtheorem{definition}[theorem]{Definition}
\newtheorem{example}[theorem]{Example}

\theoremstyle{remark}
\newtheorem{remark}[theorem]{Remark}

\numberwithin{equation}{section}

\begin{document}

\title[An extra gradient Anderson-accelerated algorithm for pseudomonotone VI\MakeLowercase{s}]{An extra gradient Anderson-accelerated algorithm for pseudomonotone variational inequalities}
\author{Xin Qu}
\address{School of Mathematics, Harbin Institute of Technology, Harbin, China}
\curraddr{Department of Applied Mathematics, The Hong Kong Polytechnic University, Hung Hom, Kowloon, Hong Kong, China}
\email{hitxin.qu@connect.polyu.hk}

\author{Wei Bian}
\address{School of Mathematics, Harbin Institute of Technology, Harbin, China}
\email{bianweilvse520@163.com}

\author{Xiaojun Chen}
\address{Department of Applied Mathematics, The Hong Kong Polytechnic University, Hung Hom, Kowloon, Hong Kong, China}
\email{xiaojun.chen@polyu.edu.hk}



\thanks{This paper is supported in part by National Natural Science Foundation of China grants (12271127, 62176073) and Hong Kong Research Grant Council grants (PolyU15300123, PolyU15300021)}
\begin{abstract}
This paper proposes an extra gradient Anderson-accelerated algorithm for solving pseudomonotone variational inequalities, which uses the extra gradient scheme with line search to guarantee the global convergence and Anderson acceleration to have fast convergent rate. We prove that the sequence generated by the proposed algorithm from any initial point converges to a solution of the pseudomonotone variational inequality problem without assuming the Lipschitz continuity and contractive condition, which are used for convergence analysis of the extra gradient method and Anderson-accelerated method, respectively in existing literatures.
Numerical experiments, particular emphasis on Harker-Pang problems, fractional programming problems, nonlinear complementarity problems, PDE problems with free boundary and linear complementarity problems, are conducted to validate the effectiveness and good performance of the proposed algorithm comparing with the extra gradient method and Anderson-accelerated method.
\end{abstract}

\maketitle


\section{Introduction}
In this paper, we consider the following variational inequality (VI) problem: find an $x^*\in\Omega$ such that
\begin{equation}\label{vi}
\langle H(x^*),x-x^*\rangle\geq0,\quad\forall x\in\Omega,
\end{equation}
where $\Omega$ is a closed convex set of $\mathbb{R}^n$ and $H:\mathbb{R}^n\rightarrow\mathbb{R}^n$ is a continuous  function, and pseudomonotone on $\Omega$, but not necessarily smooth or even Lipschitz continuous. Throughout this paper, we denote problem \eqref{vi} by ${\rm VI}(\Omega, H)$ and the solution set of \eqref{vi} by ${\rm SOL}(\Omega, H)$, and assume ${\rm SOL}(\Omega, H)\neq\emptyset.$

Variational inequalities (VIs) provide a unified framework for representing various important concepts in applied mathematics such as nonlinear equation systems, complementarity problems, optimality conditions for optimization problems and network equilibrium problems. Thus, VIs have a wide range of applications in physics, economics, engineering sciences and so on \cite{ref1,ref2,ref4,ref5,ref6}. One of the most interesting topics in VIs is to develop efficient and fast iterative algorithms to find solutions.

As a class of effective numerical methods for solving VIs, projection methods have received a lot of attention from many researchers. The earliest projection method for solving ${\rm VI}(\Omega, H)$ is the gradient projection (PG) method \cite{ref17}
\begin{equation}\label{PG}
x_{k+1}=P_{\Omega}(x_k-tH(x_k)),
\end{equation}
where $P_{\Omega}$ denotes the projection onto the set $\Omega$ and $t>0$. 
To guarantee the convergence of the PG method, it is usually assumed that either $H$ is strongly monotone and $L$-Lipschitz continuous on $\Omega$\footnote{$H$ is called $L$-Lipschitz continuous on $\Omega$ if $\|H(x)-H(y)\|\leq L\|x-y\|$ for any $x,y\in\Omega$.}, or $H$ is cocoercive.
In order to weaken the strong monotonicity of $H$, many projection methods have been proposed to solve monotone VIs \cite{new2, ref11, new1, ref10, ref13, ref37, ref16, ref14, ref15}. These methods include using some acceleration techniques to solve monotone VIs. For example, an adaptive golden ratio algorithm was proposed in \cite{new2} for solving VIs where the operator is monotone and locally Lipschitz continuous. Specifically, when the feasible set is $\mathbb{R}^n$, the monotonicity assumption of the operator in \cite{new2} can be replaced by the existence of a weak Minty solution to VIs. Under these assumptions, convergence rates were provided for both the gap function and the residual function. In \cite{new1}, the authors introduced an inertial projection and contraction algorithm and proved that, under the conditions that the operator is monotone and Lipschitz continuous, the sequence generated by the algorithm weakly converges to a solution of VIs in a Hilbert space.

Very recently, developing algorithms for different nonmonotone VIs has attracted great attention due to applications in machine learning \cite{li, ref46, ref31, ref30, ref48, ref47}. In this paper, we focus on solving ${\rm VI}(\Omega, H)$ with a pseudomonotone $H$, which is a widely used class of nonmonotone VIs. A most commonly used algorithm in the literatures for solving pseudomonotone VIs is the Korpelevich's extra gradient (EG) method \cite{ref20}.
The EG method was originally used to solve ${\rm VI}(\Omega, H)$ with a monotone and $L$-Lipschitz continuous $H$, and was later extended by Pang and Facchinei \cite{ref18} to solve the pseudomonotone VIs. After that, the EG method has been intensively studied and extended in various ways \cite{li, ref11, ref10, ref45,c1, ref49, ref46, ref22, ref37, ref31,ref30,ref48,ref15,ref23}.
It is worth pointing out that besides the $L$-Lipschitz continuity of $H$ on $\Omega$, some other conditions on $H$ are often required to guarantee the convergence of the EG method and its variants, such as Minty condition, quasimonotonicity, pseudomonotonicity and weak monotonicity.

To the best of our knowledge, current research on the Minty condition cannot simultaneously provide an analysis of both sequence convergence and convergence rates.
One of the main objectives of this paper is to design an algorithm that ensures both the sequence convergence and a fast convergence rate under the pseudomonotonicity of $H$ on $\Omega$.
In subsection \ref{section1.1}, we will review and summarize the development of EG methods for nonmonotone VIs in recent years in detail. A summary of some main comparisons is provided in {Table} \ref{ta00}.

It is known that ${\rm VI}(\Omega, H)$ is equivalent to the following fixed point problem \cite{ref18}:
\begin{equation}\label{fixed-point}
x=G(x)\coloneqq P_{\Omega}(x-tH(x)),
\end{equation}
where $t>0$. Thus, the study of PG method in \eqref{PG}, which is a fixed point method, can help to improve the performance of the algorithms for solving VIs. Anderson acceleration is efficient to improve the convergence rate of fixed point methods,  but the existing convergence analysis of Anderson acceleration requires the fixed point mapping $G$ to be contractive and piecewise smooth.  However, due to the nonmonotonicity of $H$, the mapping $G$  is not contractive, and may be even expansive, which means that we cannot directly use Anderson acceleration to solve the corresponding fixed point problem of ${\rm VI}(\Omega, H)$ in \eqref{fixed-point}. Moreover, $G$ is nonsmooth, because of  the projection operator $P_{\Omega}$.
In this paper, we introduce Anderson acceleration technique into the EG method to ensure the global sequence convergence of the algorithm and improve the convergence rate of the EG method. We will review the recent development of Anderson acceleration for fixed point problems in subsection \ref{section1.2}.
\subsection{EG methods for solving nonmonotone VIs}\label{section1.1}
The structure of the EG method proceeds as follows:
\begin{eqnarray*}
\left\{\begin{aligned}
  & x_{k+0.5}\coloneqq P_{\Omega}(x_k-tH(x_k))
  \\& x_{k+1}\coloneqq P_{\Omega}(x_k-tH(x_{k+0.5}))
\end{aligned}\right.
\end{eqnarray*}
with $t>0$.
It is important to note that the convergence results for the EG method require $H$ to be $L$-Lipschitz continuous on $\Omega$ and the stepsize $t$ to satisfy $tL<1$ \cite{ref20}. Iusem \cite{ref29} proposed a modified EG method with an updated stepsize to guarantee the efficiency of the proposed algorithm for ${\rm VI}(\Omega, H)$, in which $H$ is monotone and continuous. Recently, similar extensions have been developed not only for monotone operators but also for pseudomonotone operators \cite{li, ref31, ref30}. However, the convergence rates are not mentioned in these works. Most results on convergence rates of EG methods for VIs are established based on the $L$-Lipschitz continuity of $H$, resulting in a sublinear rate of convergence for the best-iterate of the residual term. In particular, when $H$ is monotone and $L$-Lipschitz continuous, we known that the EG method converges to a solution of ${\rm VI}(\Omega, H)$ in terms of $\min_{0\leq k\leq N}\|x_k-x_{k+0.5}\|$ with a rate of $O({1}/{\sqrt{N}})$ \cite{ref100}, which has been extended to the EG method for solving ${\rm VI}(\Omega, H)$ with a pseudomonotone and $L$-Lipschitz continuous $H$ in \cite[Lemma 12.1.10]{ref18}.

Furthermore, research on nonmonotone VIs under Minty condition has been conducted. We say ${\rm VI}(\Omega, H)$ satisfies the Minty condition if there exists an $x^*\in\Omega$ such that
\begin{equation}\label{mvi}
\langle H(x), x-x^*\rangle\geq0, \,\, \forall x\in \Omega.
\end{equation}

Under the Minty condition, Ye and He in \cite{ref47} introduced a double projection algorithm (DPA) with global convergence on the sequence, which requires computing the projection onto the intersection of a finite number of halfspaces and the closed convex set $\Omega$.
Subsequently, Lei and He in \cite{ref46} proposed a new extra gradient method (NEG) that does not involve adding halfspaces during the projection computation for solving this class of VIs under the same assumptions as in \cite{ref47}.
Then a new extra gradient type projection algorithm (NEGTP) was presented in \cite{ref48} to solve a class of continuous quasimonotone VIs satisfying $H(x)\neq\textbf{0}, \forall x\in \Omega$.
All the algorithms in \cite{ref46,ref48,ref47} have the global sequence convergence, but do not have the estimation on the convergence rate. Ye \cite{new3} proved the global convergence of the sequence for the proposed algorithm under the Minty condition. 
However, in each iteration, the algorithm in \cite{new3} requires selecting the half-space that has the largest distance from $x_k$ to some generated half-spaces. As $k$ increases, more information needs to be computed and stored.
Approximation-based Regularized Extra-gradient method (ARE), a $p^{\rm th}$-order ($p\geq1$) algorithm, was proposed in \cite{ref49} for solving monotone VIs with a convergence rate of $O(1/N^{\frac{p+1}{2}})$ on the gap function.
In \cite{zhang}, it was stated that ARE also can solve the nonmonotone VIs satisfying the Minty condition with
the convergence rate of $O(1/\sqrt{N})$ for the residual function and $O(1/N^{\frac{p}{2}})$ for the gap function. However, the algorithms in \cite{ref49,zhang} need the Lipschitz continuity of $H$ and do not have the sequence convergence on the iterates.
\begin{table}[h]
	\centering
	\label{table1}
\resizebox{\linewidth}{!}{
	\begin{tabular}{|c|c|c|c|c|c}
		
		\toprule
		  Methods  &Assumptions & Sequence convergence & Convergence rate \\
                 &&&(residual function)\\
		\midrule
		 EG \cite{ref18}    & pseudomonotone    &  $\surd$  &  $O(\frac{1}{\sqrt{N}})$ \\
                & Lipschitz continuous &&\\
		\midrule
		 ARE \cite{zhang}    & Minty condition    &  $\times$  &  $O(\frac{1}{\sqrt{N}})$ \\
              (the order $p=1$)  & Lipschitz continuous &&\\
		
        \midrule
		 DPA \cite{ref47}  & Minty condition  & $\surd$  &  $\times$ \\
		              &continuous&&\\
        \midrule
		 NEG \cite{ref46} & Minty condition  & $\surd$  &  $\times$ \\
		              &continuous&&\\
		
        \midrule
		NEGTP \cite{ref48}    & Minty condition  & $\surd$  & $\times$ \\
                      &quasimonotone&&\\
		              &continuous&&\\
                      &$H(x)\neq0, \forall x\in \Omega$&&\\

        \midrule
		EG-Anderson(1)   & pseudomonotone  & $\surd$  & $O(\frac{1}{\sqrt{N}})$ \\
		            \textbf{[This paper]}  & continuous && (locally Lipschitz continuous)\\
		\bottomrule
	\end{tabular}}
\caption{\small Summary of results on algorithms with global convergence to nonmonotone VIs
}\label{ta00}
\end{table}


\subsection{Anderson acceleration for fixed point problems}\label{section1.2}
Anderson acceleration was first proposed by Anderson in 1965 in the context of integral equations \cite{ref38}. This technique aims to improve the convergence rate of fixed point iteration by utilizing the history of search directions.
It is not necessary to compute the Jacobian of $G$, which allows it to perform effectively in various fields, including electronic structure computations \cite{ref38,ref39}, machine learning \cite{ref40}, radiation diffusion and nuclear physics \cite{ref41}.
Anderson acceleration is formally described in the following algorithm, commonly referred to as Anderson(m).


\begin{algorithm}\normalsize
\renewcommand{\thealgocf}{}
	\caption{Anderson(m)}\label{Alg}
       Choose $x_0\in \mathbb{R}^n$ and a positive integer $m$. Set $x_1=G(x_0)$ and $F_0=G(x_0)-x_0$.
      \\ \textbf{for} $k=1,2,\ldots$, \textbf{do}
       \\ set $F_k=G(x_k)-x_k$;
        \\ choose $m_k=\min\{m,k\}$;
        \\ solve
        \begin{align}\label{ass}
        \begin{split}
         \min \left\|\sum_{j=0}^{m_k} \theta_j F_{k-m_k+j}\right\| \quad \quad{\rm s.t.} \quad \sum_{j=0}^{m_k} \theta_j=1
        \end{split}
        \end{align}
       to find a solution $\{\theta_j^k:j=0,\ldots,m_k\}$, and set
        \begin{align*}
        x_{k+1}=\sum_{j=0}^{m_k}\theta_j^k G(x_{k-m_k+j}).
        \end{align*}
        \\ \textbf {end for}
\end{algorithm}
Even after a long period of use and attention, the first mathematical convergence result for Anderson acceleration had not been given until 2015 by Toth and Kelley \cite{ref42}. They showed that when $G$ is Lipschitz continuously differentiable and contractive, Anderson(m) has local r-linear convergence with r-factor $\hat{c}\in(c,1)$, and Anderson(1) has q-linear convergence with q-factor $c$, where $c$ is the contraction coefficient of the fixed point mapping. In 2019, Chen and Kelley \cite{ref39} weakened the condition of $G$, proving that this conclusion can be obtained as long as $G$ is a continuously differentiable operator. Additionally, Bian, Chen and Kelley \cite{ref43} demonstrated the q-linear convergence of Anderson(1) for general nonsmooth fixed point problems in a Hilbert space, and r-linear convergence of Anderson(m) for a special nonsmooth operator. Then, Bian and Chen \cite{addition} proved that Anderson(1) is q-linear convergent for the composite max fixed point problem with a smaller q-factor than the existing q-factors. Zhang et al. \cite{ref44} introduced a variant of Anderson acceleration that guaranteed global convergence for nonsmooth fixed point problems, but did not provide a convergence rate. Ouyang et al. \cite{newouyang}
established a globalization strategy for Anderson acceleration incorporating a nonmonotone trust-region framework. They demonstrated that the algorithm has global convergence for a class of nonexpansive mappings and showed a local r-linear convergence for contractive mappings. Moreover, the local properties of Anderson acceleration with restarting were investigated in \cite{newouyang2} in terms of function values when applied to a basic gradient scheme.

In this paper, we will utilize the Anderson method to accelerate the convergence of the EG method while ensuring the global convergence of the sequence. The contributions of this paper include the following two aspects.
 \begin{enumerate}
\item We propose a new algorithm to solve pseudomonotone ${\rm VI}(\Omega, H)$ by combining  Anderson(1) with the EG method. We prove that the sequence generated by the proposed algorithm  converges to a solution of ${\rm VI}(\Omega, H)$ without assuming the Lipschitz continuity and contractive condition of $G$.
\item Under the condition that $H$ is locally Lipschitz continuous, the convergence rate of the proposed algorithm on the residual function is not worse than the EG method. This condition is weaker than the requirement of the EG method that $H$ is Lipschitz continuous. Moreover, in numerical experiments, the proposed algorithm has been found to outperform Anderson(1) and the EG methods.
\end{enumerate}

This paper is organized as follows. In Section \ref{s2}, we briefly review some related concepts and recall some preliminary results used in this paper. In Section \ref{s3}, we use the idea of Anderson(1) to develop
an extra gradient Anderson-accelerated algorithm to solve the continuous ${\rm VI}(\Omega, H)$. Furthermore, the sequence convergence of the algorithm is analyzed and the convergence rate is provided. Finally, we use five numerical experiments to illustrate the good performance of the proposed algorithm in Section \ref{s5}.

\section{Preliminaries}\label{s2}
Let $\|\cdot\|$ denote Euclidean norm in $\mathbb{R}^n$. For a matrix $A\in\mathbb{R}^{m\times n}$, $\|A\|$ represents its $2$-norm. We begin by introducing two operators, which play a crucial role in the proposed algorithm. Additionally, we present some definitions and lemmas that will be used for the convergence analysis of the proposed algorithm.

Define the following operators $$G_t(x)\coloneqq P_{\Omega}(x-tH(x)) \quad {\rm and} \quad \tilde{G}_t(x)\coloneqq P_{\Omega}(x-tH(G_t(x))),$$
where $t>0$. Let
$$F_t(x)\coloneqq G_t(x)-x\quad \mbox{and}\quad \tilde{F}_t(x)\coloneqq \tilde{G}_t(x)-x.$$
\begin{definition}\cite{def12}\label{def}
The mapping $H:\mathbb{R}^n\rightarrow\mathbb{R}^n$ is said to be pseudomonotone on $\Omega$, if for any $x,y\in \Omega$ it holds $$\langle H(x), y-x\rangle\geq0\Rightarrow \langle H(y),y-x\rangle\geq0.$$
\end{definition}
\begin{lemma}\cite{ref}\label{lee}
For any $x\in\mathbb{R}^n$, the following statements hold.
\begin{enumerate}[{\rm (i)}]
\item $\|P_{\Omega}(x)-P_{\Omega}(y)\|^2\leq\langle P_{\Omega}(x)-P_{\Omega}(y), x-y\rangle, ~~\forall y\in\mathbb{R}^n;$

\

\item $\langle x-P_{\Omega}(x), y-P_{\Omega}(x)\rangle\leq 0, ~~\forall y\in\Omega.$
\end{enumerate}
\end{lemma}
\begin{lemma}\cite{ref45}\label{lemm1}
For any $x\in\mathbb{R}^n$ and $t_1\geq t_2>0$, the following inequalities hold:
$$\frac{\|x-P_{\Omega}(x-t_1H(x))\|}{t_1}\leq\frac{\|x-P_{\Omega}(x-t_2H(x))\|}{t_2},$$
$$\|x-P_{\Omega}(x-t_2H(x))\|\leq\|x-P_{\Omega}(x-t_1H(x))\|.$$
\end{lemma}
\begin{lemma}\cite{ref28} \label{opial}{\rm[Opial's Lemma]}
Let $S$ be a nonempty subset of $\mathbb{R}^n$, and $\{x_k\}$ a sequence of elements in $\mathbb{R}^n$. Assume that
\begin{enumerate}[{\rm (i)}]
\item every sequential cluster point of $\{x_k\}$, as $k \rightarrow \infty$, belongs to $S$;

\

\item for every $z\in S$, $\lim_{k\rightarrow\infty}\|x_k-z\|$ exists.
\end{enumerate}
Then the sequence $\{x_k\}$ converges as $k\rightarrow\infty$ to a point in $S$.
\end{lemma}
\begin{lemma}\label{lemm}{\rm\cite[Proposition 1.5.8, Exercise 1.8.29]{ref18}}
$x^*\in$ ${\rm SOL}(\Omega,H)$ if and only if it is a fixed point of $G_t$, and if and only if it is a fixed point of $\tilde{G}_t$, where $t$ can be any positive number.
\end{lemma}
\begin{lemma}\cite{c1}\label{lep}
Suppose that the mapping $H:\mathbb{R}^n\rightarrow\mathbb{R}^n$ is continuous. Then,
for all bounded sequences $\{x_k\}$, $\{y_k\}\subseteq\mathbb{R}^n$ satisfying $\lim_{k\rightarrow\infty}\|x_k-y_k\|=0$, it holds that  $\lim_{k\rightarrow\infty}\\ \|H(x_k)-H(y_k)\|=0$.
\end{lemma}
\begin{lemma}\cite{attouch}\label{attou}
Let $\{a_k\}$ and $\{\varepsilon_k\}$ be real sequences. Assume that $\{a_k\}$ is bounded from below, $\sum_{k=1}^{\infty}\varepsilon_k<\infty$ and $$a_{k+1}-a_k\leq\varepsilon_k$$ for every $k$. Then $\lim_{k\rightarrow\infty}a_k$ exists.
\end{lemma}
\section{Proposed algorithm and its convergence analysis}\label{s3}
In this section, based on Anderson acceleration and the EG methods, we propose the EG-Anderson(1) algorithm for solving ${\rm VI}(\Omega, H)$. In addition, we give the convergence analysis of this algorithm.
\subsection{Proposed algorithm}
The proposed algorithm is presented in Algorithm \ref{algright3km4}, where the line search framework in \cite{li} is used in Step 2. 
\begin{algorithm}\normalsize
    \renewcommand{\thealgocf}{1}
    \caption{EG-Anderson(1)}\label{algright3km4}
	\SetAlgoNoLine
       \textbf{Initialization:} Choose $x_0\in\Omega$, $\omega\geq0, \gamma>0$, $\tau>\frac{1}{2}$, $\rho, \mu\in(0,1)$ and $\sigma_0=1$. Give a sufficiently large $M>0$.
      \\ \textbf{for} $k=0,1,2,\ldots$, \textbf{do}
      \\ \quad \quad \textbf{Step 1:} Compute $F_{\gamma}(x_k)=P_{\Omega}(x_k-\gamma H(x_k))-x_k$.
      \\ \quad \quad \quad \quad \quad \textbf{If} $F_{\gamma}(x_k)=\textbf{0}$, then stop.
      \\ \quad \quad \quad \quad \quad \textbf{Otherwise}, let $t_k=\gamma$ and go to \textbf{Step 2}.
      \\ \quad \quad \textbf{Step 2:} Compute $y_{k+0.5} = P_{\Omega}(x_k-t_kH(x_k))$ and $y_{k+1}=P_{\Omega}(x_k-t_kH(y_{k+0.5}))$. 
       \\ \quad \quad \quad \quad \quad \textbf{If}
       \begin{equation}\label{lin}
       \begin{aligned}
       &t_k\langle H(y_{k+0.5})-H(x_k),y_{k+0.5}-y_{k+1}\rangle\\
       \leq ~&\frac{\mu}{2}\left(\|x_k-y_{k+0.5}\|^2+\|y_{k+0.5}-y_{k+1}\|^2\right),
       \end{aligned}
       \end{equation}
       \\ \quad \quad \quad \quad \quad \quad go to \textbf{Step 3}.
       \\ \quad \quad \quad \quad \quad \textbf{Otherwise}, set $t_k=\rho t_k$ and repeat \textbf{Step 2}.
       \\ \quad \quad \textbf{Step 3:} Compute $F_{t_k}(x_k)= y_{k+0.5}-x_k$ and $\tilde{F}_{t_k}(x_k)=y_{k+1}-x_k$.
       \\ \quad \quad \quad \quad \quad \textbf{If} $\|\tilde{F}_{t_k}(x_k)\|<\min\{\|F_{t_k}(x_k)\|, \omega\sigma_k^{-\tau}\}$,
       set
       \begin{equation}\label{alpha}
       \alpha_{k}= \frac{\langle \tilde{F}_{t_{k}}(x_{k}), \tilde{F}_{t_{k}}(x_{k})-F_{t_{k}}(x_{k})\rangle}{\|\tilde{F}_{t_{k}}(x_{k})-F_{t_{k}}(x_{k})\|^2},
       \end{equation}
       \\ \quad \quad \quad \quad \quad \textbf{Otherwise}, set $\alpha_k=M+1$.
       \\ \quad \quad \textbf{Step 4:} \textbf {If} $|\alpha_k|\leq M$, set
        \begin{align}\label{iteration3km}
        x_{k+1}=\alpha_k x_k+(1-\alpha_k)y_{k+1},\quad \sigma_{k+1}=\sigma_k+1.
        \end{align}
        \\ \quad \quad \quad \quad \quad \textbf{Otherwise}, set
        \begin{equation}\label{EG}
        x_{k+1}=y_{k+1},\quad \sigma_{k+1}=\sigma_k.
        \end{equation}
        \\ \textbf {end for}
\end{algorithm}

From the notations and definitions for $y_{k+0.5}$ and $y_{k+1}$ in the EG-Anderson(1) algorithm, we find that
\begin{equation*}
y_{k+0.5}=G_{t_k}(x_k)\quad\mbox{and}\quad y_{k+1}=\tilde{G}_{t_k}(x_k);
\end{equation*}
\begin{equation*}
y_{k+0.5}-x_k=F_{t_k}(x_k)\quad\mbox{and}\quad y_{k+1}-x_k=\tilde{F}_{t_k}(x_k).
\end{equation*}
\subsection{Convergence analysis}
In this subsection, we will analyze the convergence properties of the EG-Anderson(1) algorithm, including the global convergence of the sequence and the convergence rate evaluated by the residual function.
In order to categorize the iteration counts, we divide them into two subsets:
\begin{equation*}
K_{AA} = \{k_0, k_1, \cdots\}\quad\mbox{and}\quad K_{EG} = \{l_0, l_1, \cdots\},
\end{equation*}
where $K_{AA}$ consists of iterations setting by \eqref{iteration3km} and $K_{EG}$ includes the remaining iterations setting by \eqref{EG}.

If the EG-Anderson(1) algorithm is terminated in finite times, then the final output point is a solution of ${\rm VI}(\Omega, H)$. Therefore, in the following analysis we assume that the EG-Anderson(1) algorithm loops infinitely.

\begin{remark}
For $k_i\in K_{AA}$, note that $\|\tilde{F}_{t_{k_i}}(x_{k_i})-F_{t_{k_i}}(x_{k_i})\|\neq0$ due to $\|\tilde{F}_{t_{k_i}}(x_{k_i})\|<\|F_{t_{k_i}}(x_{k_i})\|,$ thus $\alpha_{k_i}$ is well-defined. Moreover, $\alpha_{k_i}$ is the optimal solution of $$\min \left\|\alpha F_{t_{k_i}}(x_{k_i})+(1-\alpha)\tilde{F}_{t_{k_i}}(x_{k_i})\right\|.$$
\end{remark}
We start the convergence analysis of the EG-Anderson(1) algorithm by proving that \eqref{lin} terminates after a finite number of loops.
\begin{lemma}
The EG-Anderson(1) algorithm is well-defined.
\end{lemma}
{\it Proof}\quad We will show that the EG-Anderson(1) algorithm is well-defined by proving that for every $k$ there exists $t_k$ satisfying \eqref{lin} when $x_k\notin {\rm SOL}(\Omega, H)$.

From the updated form of $t_k$ in the EG-Anderson(1) algorithm, it can be reformulated as $t_k=\gamma \rho^{m_k}$, where $m_k$ is the smallest nonnegative integer $m$ satisfying
 \begin{equation}\label{add}
       \begin{aligned}
       &\gamma\rho^m\left \langle H\left(y^{(m)}_{k+0.5}\right)-H\left(x_k\right),y^{(m)}_{k+0.5}-y^{(m)}_{k+1}\right\rangle\\
       \leq ~ &\frac{\mu}{2}\left(\left\|x_k-y^{(m)}_{k+0.5}\right\|^2+\left\|y^{(m)}_{k+0.5}-y^{(m)}_{k+1}\right\|^2\right),
       \end{aligned}
       \end{equation}
where $y^{(m)}_{k+0.5}\coloneqq P_{\Omega}\left(x_k-\gamma \rho^{m}H\left(x_k\right)\right)$ and $y^{(m)}_{k+1}\coloneqq P_{\Omega}\left(x_k-\gamma \rho^{m}H\left(y^{(m)}_{k+0.5}\right)\right).$

	If there exists a nonnegative integer $\bar{m}$ such that $y^{(\bar{m})}_{k+0.5} = y^{(\bar{m})}_{k+1}$, then there must exist an integer $m_k$ between 0 and $\bar{m}$ such that \eqref{add} holds.
We consider the situation $y^{(m)}_{k+0.5} \neq y^{(m)}_{k+1}$ for any nonnegative integer $m$ and assume the contrary that for all $m$ we have
\begin{equation}\label{zhu1}
\gamma\rho^m\left\langle H\left(y^{(m)}_{k+0.5}\right)-H(x_k),y^{(m)}_{k+0.5}-y^{(m)}_{k+1}\right\rangle >\frac{\mu}{2}\left(\left\|x_k-y^{(m)}_{k+0.5}\right\|^2+\left\|y^{(m)}_{k+0.5}-y^{(m)}_{k+1}\right\|^2\right).
\end{equation}
On one hand, by Cauchy–Schwarz inequality, we obtain
\begin{equation}\label{dw12}
\gamma\rho^{m}\left\langle H\left(y^{(m)}_{k+0.5}\right)-H\left(x_k\right), y^{(m)}_{k+0.5}-y^{(m)}_{k+1}\right\rangle\leq \gamma\rho^{m}\left\|H\left(y^{(m)}_{k+0.5}\right)-H(x_k)\right\|\left\|y^{(m)}_{k+0.5}-y^{(m)}_{k+1}\right\|.
\end{equation}
On the other hand, we also find
\begin{equation}\label{dw22}
\left\|x_k-y^{(m)}_{k+0.5}\right\|^2+\left\|y^{(m)}_{k+0.5}-y^{(m)}_{k+1}\right\|^2\geq 2\left\|x_k-y^{(m)}_{k+0.5}\right\|\left\|y^{(m)}_{k+0.5}-y^{(m)}_{k+1}\right\|.
\end{equation}
Combining \eqref{zhu1} with \eqref{dw12} and \eqref{dw22}, we deduce that
\begin{equation}\label{eq222}
\frac{\left\|x_k-y^{(m)}_{k+0.5}\right\|}{\gamma \rho^m}\leq \frac{1}{\mu}\left\|H\left(y^{(m)}_{k+0.5}\right)-H(x_k)\right\|.
\end{equation}
Since $x_k\notin {\rm SOL}(\Omega, H)$, we discuss the following two cases.

(i) If $x_k\in\Omega$, from the definition of $y^{(m)}_{k+0.5}$ and the continuity of $P_{\Omega}$, we have $$\lim_{m\rightarrow\infty}\left\|x_k-y^{(m)}_{k+0.5}\right\|=0.$$ In view of the continuity of $H$, we get $\lim_{m\rightarrow\infty}\left\|H(x_k)-H\left(y^{(m)}_{k+0.5}\right)\right\|=0$. This together with \eqref{eq222} yields
\begin{equation}\label{d3}
\lim_{m\rightarrow\infty}\frac{\left\|x_k-y^{(m)}_{k+0.5}\right\|}{\gamma \rho^m}=0.
\end{equation}
By the definition of $y^{(m)}_{k+0.5}$ and using Lemma \ref{lee}-(ii), we get $$\left\langle y^{(m)}_{k+0.5}-x_k+\gamma \rho^m H(x_k), x-y^{(m)}_{k+0.5}\right\rangle\geq0, \forall x\in\Omega,$$
which implies
\begin{equation}\label{d2}
\left \langle\frac{y^{(m)}_{k+0.5}-x_k}{\gamma \rho^m}+H(x_k), x-y^{(m)}_{k+0.5}\right \rangle\geq0, \forall x\in\Omega.
\end{equation}
Taking the limit $m\rightarrow\infty$ in \eqref{d2} and using \eqref{d3} and $\lim_{m\rightarrow\infty}y^{(m)}_{k+0.5}=x_k$, we obtain $\langle H(x_k), x-x_k\rangle\geq0, \forall x\in\Omega.$ It can be deduced that $x_k\in {\rm SOL}(\Omega, H)$ and this leads to a contraction.

(ii) If $x_k\notin\Omega,$ we can conclude that
\begin{equation}\label{o1}
\lim_{m\rightarrow\infty}\left\|x_k-y^{(m)}_{k+0.5}\right\|=\|x_k-P_{\Omega}(x_k)\|>0,
\end{equation}
and
\begin{equation}\label{o2}
\lim_{m\rightarrow\infty}\gamma\rho^m\left\|H\left(y^{(m)}_{k+0.5}\right)-H(x_k)\right\|=0.
\end{equation}
Rearranging the terms in \eqref{eq222}, we find $$\left\|x_k-y^{(m)}_{k+0.5}\right\|\leq\frac{1}{\mu}\gamma\rho^m\left(\left\|H\left(y^{(m)}_{k+0.5}\right)-H(x_k)\right\|\right).$$
\\Taking the limit $m\rightarrow\infty$ in the above inequality, we can find a contradiction with \eqref{o1} and \eqref{o2}. Hence, the proof is fully established.
\qed

The next two lemmas are instrumental in establishing the key findings of this section.
\begin{lemma}\label{oa}
Let $\{x_k\}$ be the sequence generated by the EG-Anderson(1) algorithm. Then we have
$$\left(1-\sqrt{\frac{\mu}{2-\mu}}\right)\|F_{t_k}(x_k)\|\leq\|\tilde{F}_{t_k}(x_k)\|\leq \left(1+\sqrt{\frac{\mu}{2-\mu}}\right)\|F_{t_k}(x_k)\|.$$
\end{lemma}
{\it Proof} \quad For any $k$, by the condition of $t_k$ in \eqref{lin} and Lemma \ref{lee}-(i), we obtain
\begin{equation*}
\begin{split}
\|y_{k+0.5}-y_{k+1}\|^2 =~ &\|P_{\Omega}(x_k-t_kH(x_k))-P_{\Omega}(x_k-t_kH(y_{k+0.5}))\|^2
\\ \leq~ & \langle y_{k+0.5}-y_{k+1}, t_k\left(H(y_{k+0.5})-H(x_k)\right)\rangle
\\ \leq~ & \frac{\mu}{2}\|x_k-y_{k+0.5}\|^2+\frac{\mu}{2}\|y_{k+0.5}-y_{k+1}\|^2,
\end{split}
\end{equation*}
which implies
\begin{equation}\label{pm}
\|y_{k+1}-y_{k+0.5}\|^2\leq\frac{\mu}{2-\mu}\|x_k-y_{k+0.5}\|^2.
\end{equation}
Since $\mu\in(0,1)$, then $\frac{\mu}{2-\mu}\in(0,1)$. From the triangle inequality and \eqref{pm}, we have
\begin{equation*}
\begin{split}
\|\tilde{F}_{t_k}(x_k)\|=~&\|y_{k+1}-x_k\|\geq\|y_{k+0.5}-x_k\|-\|y_{k+1}-y_{k+0.5}\|
\\ \geq~ & \|y_{k+0.5}-x_k\|-\sqrt{\frac{\mu}{2-\mu}}\|y_{k+0.5}-x_k\|
\\ = ~& \left(1-\sqrt{\frac{\mu}{2-\mu}}\right)\|F_{t_k}(x_k)\|
\end{split}
\end{equation*}
and
\begin{equation*}
\begin{split}
\|\tilde{F}_{t_k}(x_k)\|=~&\|y_{k+1}-x_k\|\leq\|y_{k+0.5}-x_k\|+\|y_{k+1}-y_{k+0.5}\|
\\ \leq ~& \left(1+\sqrt{\frac{\mu}{2-\mu}}\right)\|F_{t_k}(x_k)\|.
\end{split}
\end{equation*}
The proof is completed.
\qed
\begin{lemma}\label{lemm3}
Let $\{x_k\}$ be the sequence generated by the EG-Anderson(1) algorithm and $x^*\in{\rm SOL}(\Omega, H)$. For every $k$, it holds that
\begin{equation*}
\|y_{k+1}-x^*\|^2\leq\|x_{k}-x^*\|^2-(1-\mu)\|y_{k+0.5}-x_{k}\|^2-(1-\mu)\|y_{k+1}-y_{k+0.5}\|^2.
\end{equation*}
\end{lemma}
{\it Proof} \quad In view of the pseudomonotonicity of $H$ and $x^*\in {\rm SOL}(\Omega, H)$, we deduce that $\langle t_kH(y_{k+0.5}), x^*-y_{k+0.5}\rangle\leq0$, which gives
\begin{equation}\label{bh}
\begin{split}
\langle t_kH(y_{k+0.5}), x^*-y_{k+1}\rangle\leq \langle t_kH(y_{k+0.5}), y_{k+0.5}-y_{k+1}\rangle.
\end{split}
\end{equation}
By the definition of $y_{k+0.5}$, Lemma \ref{lee}-(ii) and \eqref{lin}, we obtain
\begin{equation}\label{bh1}
\begin{split}
& \langle x_k-t_kH(y_{k+0.5})-y_{k+0.5}, y_{k+1}-y_{k+0.5}\rangle
\\ =~& \langle x_k-t_kH(x_k)-y_{k+0.5}, y_{k+1}-y_{k+0.5}\rangle + t_k\langle H(x_k)-H(y_{k+0.5}), y_{k+1}-y_{k+0.5}\rangle
\\\overset{\eqref{lin}}{\leq} & \langle x_k-t_kH(x_k)-P_{\Omega}(x_k-t_kH(x_k)), y_{k+1}-P_{\Omega}(x_k-t_kH(x_k))\rangle
\\ & + \frac{\mu}{2}\|x_k-y_{k+0.5}\|^2+\frac{\mu}{2}\|y_{k+0.5}-y_{k+1}\|^2
\\ \leq~ & \frac{\mu}{2}\|x_k-y_{k+0.5}\|^2+\frac{\mu}{2}\|y_{k+0.5}-y_{k+1}\|^2.
\end{split}
\end{equation}
Based on the definition of $y_{k+1}$, Lemma \ref{lee}-(ii), \eqref{bh} and \eqref{bh1}, and using the same idea as in {\rm\cite[Lemma 3.3 ]{li}}, we conclude that
\begin{equation*}
\begin{split}
\|y_{k+1}-x^*\|^2
 \leq~ & \|x_k-x^*\|^2-\|x_k-y_{k+1}\|^2+2\langle t_kH(y_{k+0.5}), x^*-y_{k+1}\rangle
\\ \overset{\eqref{bh}}{\leq} & \|x_k-x^*\|^2-\|x_k-y_{k+1}\|^2+2\langle t_kH(y_{k+0.5}), y_{k+0.5}-y_{k+1}\rangle
\\ = ~ & \|x_k-x^*\|^2-\|x_k-y_{k+0.5}\|^2-\|y_{k+0.5}-y_{k+1}\|^2
\\ & +2\langle x_k-t_kH(y_{k+0.5})-y_{k+0.5}, y_{k+1}-y_{k+0.5}\rangle
\\ \overset{\eqref{bh1}}{\leq} & \|x_k-x^*\|^2-(1-\mu)\|x_k-y_{k+0.5}\|^2-(1-\mu)\|y_{k+0.5}-y_{k+1}\|^2.
\end{split}
\end{equation*}
\qed

Now, utilizing the Opial's Lemma, we can state and prove our main convergence result in what follows.
\begin{theorem}\label{lem3}
Let $\{x_k\}$ be the sequence generated by the EG-Anderson(1) algorithm. Then the sequence $\{x_k\}$ converges to a solution of ${\rm VI}(\Omega, H)$.
\end{theorem}
{\it Proof} \quad We will prove this theorem from three steps.

Step 1: $\{x_k\}$ is bounded.

Let $x^*$ be a solution of problem ${\rm VI}(\Omega, H)$. If $k_i\in K_{AA}$, let $\beta_{k_i}\coloneqq 1-\alpha_{k_i}$ and by the definition of $x_{k_i+1}$ in \eqref{iteration3km}, we know that
\begin{equation}\label{kun1111}
\begin{split}
\|x_{k_i+1}-x^*\|^2 = ~& \|\alpha_{k_i} x_{k_i}+\beta_{k_i} y_{k_i+1}-x^*\|^2
\\ = ~& \alpha_{k_i}^2\|x_{k_i}-x^*\|^2+\beta_{k_i}^2\|y_{k_i+1}-x^*\|^2+2\alpha_{k_i}\beta_{k_i}\langle x_{k_i}-x^*, y_{k_i+1}-x^*\rangle
\\  = ~& \alpha_{k_i}^2\|x_{k_i}-x^*\|^2+\beta_{k_i}^2\|y_{k_i+1}-x^*\|^2
\\ &+2\alpha_{k_i}\beta_{k_i} \left(\frac{1}{2}\|x_{k_i}-x^*\|^2+\frac{1}{2}\|y_{k_i+1}-x^*\|^2-\frac{1}{2}\|x_{k_i}-y_{k_i+1}\|^2\right)
\\ =~ & \alpha_{k_i}\|x_{k_i}-x^*\|^2+\beta_{k_i}\|y_{k_i+1}-x^*\|^2-\alpha_{k_i}\beta_{k_i}\|x_{k_i}-y_{k_i+1}\|^2.
\end{split}
\end{equation}
From Lemma \ref{lemm3}, we can obtain
\begin{equation}\label{haha}
\|y_{k_i+1}-x^*\|^2\leq\|x_{k_i}-x^*\|^2-(1-\mu)\|y_{k_i+0.5}-x_{k_i}\|^2-(1-\mu)\|y_{k_i+1}-y_{k_i+0.5}\|^2.
\end{equation}
Since $\|\tilde{F}_{t_{k_i}}(x_{k_i})\|<\|F_{t_{k_i}}(x_{k_i})\|$ when $k_i\in K_{AA}$, then $$\langle F_{t_{k_i}}(x_{k_i}), \tilde{F}_{t_{k_i}}(x_{k_i})\rangle<\|F_{t_{k_i}}(x_{k_i})\|^2,$$ which implies
$$\langle \tilde{F}_{t_{k_i}}(x_{k_i}), \tilde{F}_{t_{k_i}}(x_{k_i})-F_{t_{k_i}}(x_{k_i})\rangle<
\|\tilde{F}_{t_{k_i}}(x_{k_i})-F_{t_{k_i}}(x_{k_i})\|^2.$$
By \eqref{alpha}, we have
$\alpha_{k_i}<1.$ Thus $\beta_{k_i}=1-\alpha_{k_i}>0$.

Introducing \eqref{haha} into \eqref{kun1111}, we deduce that
\begin{equation}\label{e8}
\begin{split}
\|x_{k_i+1}-x^*\|^2\leq ~& \|x_{k_i}-x^*\|^2 -(1-\mu)\beta_{k_i}\|y_{k_i+0.5}-x_{k_i}\|^2
\\ &-(1-\mu)\beta_{k_i}\|y_{k_i+0.5}-y_{k_i+1}\|^2-\alpha_{k_i}\beta_{k_i}\|x_{k_i}-y_{k_i+1}\|^2.
\end{split}
\end{equation}
In view of $|\alpha_{k_i}|\leq M$ and $\beta_{k_i}=1-\alpha_{k_i}$, we get $|\beta_{k_i}|\leq M+1$. This together with \eqref{e8}, $\mu\in(0,1)$, $\beta_{k_i}>0$ and $\|x_{k_i}-y_{k_i+1}\|=\|\tilde{F}_{t_{k_i}}(x_{k_i})\|\leq \omega \sigma_{k_i}^{-\tau}=\omega (\sigma_0+i)^{-\tau}=\omega (1+i)^{-\tau}$, we deduce that
\begin{equation}\label{qw8}
\begin{split}
\|x_{k_i+1}-x^*\|^2 & \leq \|x_{k_i}-x^*\|^2+|\alpha_{k_i}||\beta_{k_i}|\|\tilde{F}_{t_{k_i}}(x_{k_i})\|^2
\\ & \leq \|x_{k_i}-x^*\|^2 + \epsilon_{k_i},
\end{split}
\end{equation}
where $\epsilon_{k_i}\coloneqq M(M+1)\omega^2(1+i)^{-2\tau}.$

If $l_j\in K_{EG}$, by Lemma \ref{lemm3} and $\mu\in(0,1)$, we conclude that
\begin{equation}\label{haha11}
\begin{split}
 \|x_{l_j+1}-x^*\|^2=~& \|y_{l_j+1}-x^*\|^2
\\ \leq~ & \|x_{l_j}-x^*\|^2-(1-\mu)\|y_{l_j+0.5}-x_{l_j}\|^2-(1-\mu)\|y_{l_j+1}-y_{l_j+0.5}\|^2
\\  \leq ~& \|x_{l_j}-x^*\|^2.
\end{split}
\end{equation}
By defining $\epsilon_{l_j}=0$ and combining \eqref{qw8} and \eqref{haha11}, for every $k$, we find
\begin{equation}\label{a1}
\|x_{k+1}-x^*\|^2\leq\|x_{k}-x^*\|^2 + \epsilon_k,
\end{equation}
where $\sum_{k=0}^\infty \epsilon_k=\sum_{i=0}^\infty \epsilon_{k_i}<\infty$ because of $\tau>\frac{1}{2}$. 
 From \eqref{a1}, we further have
\begin{equation*}\label{afg}
\|x_k-x^*\|^2\leq\|x_0-x^*\|^2 + \sum_{k=0}^\infty \epsilon_k <\infty.
\end{equation*}
Hence $\{x_{k}\}$ is bounded.

Step 2: any cluster point of $\{x_k\}$ belongs to ${\rm SOL}(\Omega,H)$, i.e. a solution of ${\rm VI}(\Omega, H)$.

By the boundedness of $\{x_{k}\}$, it has at least one cluster point, denoted by $\bar{x}$ with convergence subsequence $\{x_{p_{j}}\}$ satisfying $\lim_{j\rightarrow\infty} x_{p_{j}}=\bar{x}$.
Next, we will prove $\bar{x}\in{\rm SOL}(\Omega,H)$ from the following two steps.

Step 2.1:
 $\lim_{k\rightarrow\infty}\|F_{t_k}(x_k)\|=0.$

We know that $K_{AA}\cup K_{EG}$ is infinite. The following proof is divided into two cases for discussion.
\\ (i) We first consider the case that both $K_{AA}$ and $K_{EG}$ are infinite.
Rearranging the terms in \eqref{haha11} and using $\mu\in(0,1)$, we infer that
\begin{equation}\label{hak}
(1-\mu)\|F_{t_{l_j}}(x_{l_j})\|^2=(1-\mu)\|y_{l_j+0.5}-x_{l_j}\|^2\leq\|x_{l_j}-x^*\|^2-\|x_{l_j+1}-x^*\|^2.
\end{equation}
Combining \eqref{qw8} and \eqref{hak}, we deduce that
\begin{equation}\label{slk}
\sum_{j=0}^{\infty}\|F_{t_{l_j}}(x_{l_j})\|^2\leq\frac{1}{1-\mu}\|x_0-x^*\|^2+\frac{1}{1-\mu}\sum_{i=0}^{\infty}\epsilon_{k_i}<\infty.
\end{equation}
Moreover, we know $\|\tilde{F}_{t_{k_i}}(x_{k_i})\|\leq \omega(1+i)^{-\tau}$ when $k_i\in K_{AA}$. Together with Lemma \ref{oa}, we have $$\|F_{t_{k_i}}(x_{k_i})\|^2\leq\left(1-\sqrt{\frac{\mu}{2-\mu}}\right)^{-2}\|\tilde{F}_{t_{k_i}}(x_{k_i})\|^2\leq \left(1-\sqrt{\frac{\mu}{2-\mu}}\right)^{-2}\omega^2(1+i)^{-2\tau}.$$ Thus
\begin{equation}\label{cx}
\sum_{i=0}^{\infty}\|F_{t_{k_i}}(x_{k_i})\|^2\leq \left(1-\sqrt{\frac{\mu}{2-\mu}}\right)^{-2}\omega^2\sum_{i=0}^{\infty}(1+i)^{-2\tau}<\infty.
\end{equation}
Since $\tau>\frac{1}{2}$, there exists a $C_0>0$ such that $\sum_{i=0}^{\infty}(1+i)^{-2\tau}\leq C_0.$
Combining \eqref{slk} and \eqref{cx}, we deduce that
\begin{equation}\label{eq6}
\sum_{k=0}^{\infty}\|F_{t_k}(x_k)\|^2=\sum_{i=0}^{\infty}\|F_{t_{k_i}}(x_{k_i})\|^2+\sum_{j=0}^{\infty}\|F_{t_{l_j}}(x_{l_j})\|^2\leq C,
\end{equation}
where $C\coloneqq \left(\left(1-\sqrt{\frac{\mu}{2-\mu}}\right)^{-2}+\frac{1}{1-\mu}M(M+1)\right)\omega^2C_0+\frac{1}{1-\mu}\|x_0-x^*\|^2.$ Hence we conclude that $\lim_{k\rightarrow\infty}\|F_{t_k}(x_{k})\|=0$. 
\\(ii) 
	We then consider the case that either $K_{AA}$ or $K_{EG}$ is finite. 
 In this situation, the proof becomes simpler, as we only need to use either \eqref{slk} or \eqref{cx} to obtain $\sum_{k=0}^{\infty}\|F_{t_k}(x_k)\|^2<\infty$.  Therefore, we have $\lim_{k\rightarrow\infty}\|F_{t_k}(x_{k})\|=0$.


Step 2.2:
$\lim_{k\rightarrow\infty}\|F_{\gamma}(x_k)\|=0.$

	

Let $\tilde{y}_{k+0.5}\coloneqq P_{\Omega}(x_{k}-t_{k}\rho^{-1}H(x_{k}))$. Applying Lemma \ref{lemm1} and $t_{k}\rho^{-1}>t_{k},$ we get
\begin{equation}\label{eq1}
\|x_{k}-\tilde{y}_{k+0.5}\|\leq \rho^{-1}\|x_{k}-y_{k+0.5}\|=\rho^{-1}\|F_{t_{k}}(x_{k})\|,
\end{equation}
which implies $\lim_{k\rightarrow\infty}\|x_{k}-\tilde{y}_{k+0.5}\|=0.$

Below, we will estimate $\frac{1}{t_{k}}\|x_{k}-\tilde{y}_{k+0.5}\|$ by dividing it into two cases based on the value of $t_k$.
\\ (i) If $t_k=\gamma,$ we have 
\begin{equation}\label{review}
\frac{1}{t_{k}}\|x_{k}-\tilde{y}_{k+0.5}\|=\frac{1}{\gamma}\|x_{k}-\tilde{y}_{k+0.5}\|.
\end{equation}	
\\(ii) If $t_k<\gamma$,
let $\tilde{y}_{k+1}\coloneqq P_{\Omega}(x_{k}-t_{k}\rho^{-1}H(\tilde{y}_{k+0.5}))$.
From the condition of $t_{k}$ in \eqref{lin}, we know that $t_{k}\rho^{-1}$ satisfies
\begin{equation}\label{zhu}
t_{k}\rho^{-1}\left\langle H(\tilde{y}_{k+0.5})-H(x_{k}), \tilde{y}_{k+0.5}-\tilde{y}_{k+1}\right\rangle>\frac{\mu}{2}\left(\left\|x_{k}-\tilde{y}_{k+0.5}\right\|^2+\left\|\tilde{y}_{k+0.5}-\tilde{y}_{k+1}\right\|^2\right),
\end{equation}
which implies $\tilde{y}_{k+0.5}\neq \tilde{y}_{k+1}.$
Then for \eqref{zhu}, based on a similar estimation as for \eqref{zhu1}, we can conclude that
\begin{equation}\label{eq2}
\frac{1}{t_{k}}\|x_{k}-\tilde{y}_{k+0.5}\|<\mu^{-1}\rho^{-1}\|H(x_{k})-H(\tilde{y}_{k+0.5})\|.
\end{equation}
Combining \eqref{review} and \eqref{eq2}, we can obtain
\begin{equation}\label{hap}
\frac{1}{t_{k}}\|x_{k}-\tilde{y}_{k+0.5}\|\leq\max\left\{\frac{1}{\gamma}\|x_{k}-\tilde{y}_{k+0.5}\|, \mu^{-1}\rho^{-1}\|H(x_{k})-H(\tilde{y}_{k+0.5})\|\right\}.
\end{equation}
Since $\{x_{k}\}$ is bounded and $\lim_{k\rightarrow\infty}\|x_{k}-\tilde{y}_{k+0.5}\|=0$, we find $\{\tilde{y}_{k+0.5}\}$ is bounded. Combining the continuity of $H$ and Lemma \ref{lep}, we conclude that $$\lim_{k\rightarrow\infty}\|H(x_{k})-H(\tilde{y}_{k+0.5})\|=0.$$ This together with $\lim_{k\rightarrow\infty}\|x_{k}-\tilde{y}_{k+0.5}\|=0$ and \eqref{hap} yields
$$\lim_{k\rightarrow\infty}\frac{1}{t_{k}}\|x_{k}-\tilde{y}_{k+0.5}\|=0.$$
Again by Lemma \ref{lemm1} and $t_{k}\rho^{-1}>t_{k},$ we get
\begin{equation}\label{eq3}
\|x_{k}-\tilde{y}_{k+0.5}\|\geq\|x_{k}-y_{k+0.5}\|.
\end{equation}
Then, we obtain$$\lim_{k\rightarrow\infty}\frac{1}{t_{k}}\|x_{k}-y_{k+0.5}\|\leq\lim_{k\rightarrow\infty}\frac{1}{t_{k}}\|x_{k}-\tilde{y}_{k+0.5}\|=0,$$
which means $\lim_{k\rightarrow\infty}\frac{1}{t_{k}}\|F_{t_{k}}(x_{k})\|=\lim_{k\rightarrow\infty}\frac{1}{t_{k}}\|x_{k}-y_{k+0.5}\|=0.$ Hence, we conclude that $\lim_{k\rightarrow\infty}\frac{1}{t_k}\|F_{t_k}(x_{k})\|=0$. 


Recalling Lemma \ref{lemm1} and $t_k\leq\gamma$, we have $\frac{1}{t_k}\|F_{t_k}(x_k)\|\geq\frac{1}{\gamma}\|F_{\gamma}(x_k)\|$ and $\|F_{t_k}(x_k)\|\leq\|F_{\gamma}(x_k)\|.$ Then we conclude that
\begin{equation}\label{eq4}
\|F_{t_k}(x_k)\|\leq\|F_{\gamma}(x_k)\|\leq\frac{\gamma}{t_k}\|F_{t_k}(x_k)\|,
\end{equation}
which gives $$\lim_{k\rightarrow\infty}\|F_{\gamma}(x_k)\|=0.$$
Together with $\lim_{j\rightarrow\infty} x_{p_{j}}= \bar{x}$ and Lemma \ref{lemm}, we deduce that $\bar{x}\in{\rm SOL}(\Omega,H)$.

Step 3: $\{x_k\}$ is convergent to a solution of ${\rm VI}(\Omega, H)$.

Since $\|x_k-x^*\|^2$ is nonnegative, $\sum_{k=0}^{\infty}\epsilon_k<\infty$ and \eqref{a1} holds, by applying Lemma \ref{attou} with $a_k=\|x_k-x^*\|^2$, $\lim_{k\rightarrow\infty}\|x_k-x^*\|^2$ exists for any $x^*\in{\rm SOL}(\Omega,H)$. Together this with Step 2, the proof is completed by Lemma \ref{opial}.
\qed

If $H$ is also locally Lipschitz continuous at any solution of ${\rm VI}(\Omega, H)$, we derive the following conclusion about the convergence rate on the residual function.
\begin{theorem}(Best-iterate convergence rate)
Suppose that $H$ is locally Lipschitz continuous at any solution of ${\rm VI}(\Omega, H)$. Let $\{x_k\}$ be the sequence generated by the EG-Anderson(1) algorithm. Then there exists a positive integer $N_0$ such that
$$\min_{N_0+1\leq k\leq N}\|F_{\gamma}(x_k)\|^2=O\left(\frac{1}{N}\right).$$
\end{theorem}
{\it Proof}\quad 
	By Theorem \ref{lem3}, the sequence $\{x_k\}$ converges to a solution $x^*$ of ${\rm VI}(\Omega, H)$ and $\lim_{k\rightarrow\infty}\|x_{k}-\tilde{y}_{k+0.5}\|=0.$ Thus, we know that the sequence $\{\tilde{y}_{k+0.5}\}$ also converges to $x^*$. From the locally Lipschitz continuity of $H$, there exist an $r\in(0,1)$, an $L^*>0$ and a positive integer $N_0$ such that for $k>N_0$, we have $\|x_{k}-x^*\|\leq r$ and $$\|H(x_{k})-H(\tilde{y}_{k+0.5})\|\leq L^*\|x_{k}-\tilde{y}_{k+0.5}\|.$$
 First, we will discuss the relationship between $\|F_{\gamma}(x_k)\|$ and $\|F_{t_{k}}(x_{k})\|$ in two cases.
\\(i) For $t_k=\gamma,$ we have
\begin{equation}\label{m1}
\|F_{\gamma}(x_{k})\|=\|F_{t_{k}}(x_{k})\|.
\end{equation}
(ii) For $t_k<\gamma,$
combining \eqref{eq1}, \eqref{eq2}, \eqref{eq3} and \eqref{eq4} yields that
\begin{equation}\label{eq5}
\begin{split}
\|F_{\gamma}(x_{k})\| & \overset{\eqref{eq4}}\leq\frac{\gamma}{t_{k}}\|F_{t_{k}}(x_{k})\|
\\ & \overset{\eqref{eq3}}\leq\frac{\gamma}{t_{k}}\|x_{k}-\tilde{y}_{{k}+0.5}\|
\\ & \overset{\eqref{eq2}}<\gamma\rho^{-1}\mu^{-1}\|H(x_{k})-H(\tilde{y}_{{k}+0.5})\|
\\ & \hspace{0.6em}\leq~\gamma\rho^{-1}\mu^{-1}L^*\|x_{k}-\tilde{y}_{{k}+0.5}\|
\\ & \overset{\eqref{eq1}}\leq\gamma\rho^{-2}\mu^{-1}L^*\|F_{t_{k}}(x_{k})\|, \quad \forall k>N_0.
\end{split}
\end{equation}
Combining \eqref{m1} and \eqref{eq5}, we know that for any $k>N_0$, we have
\begin{equation}\label{eqq1}
 \|F_{\gamma}(x_{k})\|\leq\max\{\gamma\rho^{-2}\mu^{-1}L^*,1\}\|F_{t_{k}}(x_{k})\|.
 \end{equation}
Next, we will prove that $\sum_{k=N_0+1}^{\infty}\|F_{t_k}(x_k)\|^2$ is finite. 
 We only consider the case that both $K_{AA}$ and $K_{EG}$ are infinite, as the analysis is similar or simpler when either of them is finite.
For the aforementioned $N_0$, there exist $i_0$ and $j_0$ such that $\{k_i:i\geq i_0\}\cup\{l_j:j\geq j_0\}=\{k: k\geq N_0+1\}.$ Then we know
\begin{equation}\label{pig}
\sum_{k=N_0+1}^{\infty}\|F_{t_k}(x_k)\|^2 = \sum_{i=i_0}^{\infty}\|F_{t_{k_i}}(x_{k_i})\|^2+\sum_{j=j_0}^{\infty}\|F_{t_{l_j}}(x_{l_j})\|^2.
\end{equation}

From \eqref{qw8} and \eqref{hak}, we get
\begin{equation}\label{dog}
0\leq\sum_{i=i_0}^{\infty}\left(\|x_{k_i}-x^*\|^2-\|x_{k_i+1}-x^*\|^2\right)+\sum_{i=i_0}^{\infty}\epsilon_{k_i}
\end{equation}
and
\begin{equation}\label{dog2}
(1-\mu)\sum_{j=j_0}^{\infty}\|F_{t_{l_j}}(x_{l_j})\|^2\leq\sum_{j=j_0}^{\infty}\left(\|x_{l_j}-x^*\|^2-\|x_{l_j+1}-x^*\|^2\right),
\end{equation}
respectively.

Adding \eqref{dog} and \eqref{dog2}, and using $\|x_{N_0+1}-x^*\|\leq r<1$, we obtain
\begin{equation}\label{slk9}
\begin{split}
\sum_{j=j_0}^{\infty}\|F_{t_{l_j}}(x_{l_j})\|^2\leq & \frac{1}{1-\mu}\|x_{N_0+1}-x^*\|^2+\frac{1}{1-\mu}\sum_{i=0}^{\infty}\epsilon_{k_i}
\\ \leq & \frac{1}{1-\mu} +\frac{1}{1-\mu}M(M+1)\omega^2C_0.
\end{split}
\end{equation}
By \eqref{cx}, \eqref{pig} and \eqref{slk9}, we conclude that
\begin{equation*}
\sum_{k=N_0+1}^{\infty}\|F_{t_k}(x_k)\|^2 \leq \tilde{C},
\end{equation*}
where $\tilde{C}:= \left(\left(1-\sqrt{\frac{\mu}{2-\mu}}\right)^{-2}+\frac{1}{1-\mu}M(M+1)\right)\omega^2C_0+\frac{1}{1-\mu}.$

This together with \eqref{eqq1} yields
\begin{equation*}\label{eq62}
\begin{split}
\sum_{k=N_0+1}^{\infty}\|F_{\gamma}(x_k)\|^2 & \leq\max\{\gamma\rho^{-2}\mu^{-1}L^*,1\}^2\sum_{k=N_0+1}^{\infty}\|F_{t_k}(x_k)\|^2\leq C^*
\end{split}
\end{equation*}
with $C^*\coloneqq \max\{\gamma\rho^{-2}\mu^{-1}L^*,1\}^2\tilde{C}.$ Hence we conclude that $$\min_{N_0+1\leq k\leq N}\|F_{\gamma}(x_k)\|^2\leq\frac{1}{N-N_0}\sum_{k=N_0+1}^N\|F_{\gamma}(x_k)\|^2\leq\frac{1}{N-N_0}C^*.$$
\qed
\begin{remark}
If the pseudomonotone operator $H$ is $L$-Lipschitz continuous, we will no longer need the line search step \eqref{lin} in the EG-Anderson(1) algorithm, in which case we can take $t_k$ to be the constant $t>0$ satisfying $tL<1$. At this situation, let $\{x_k\}$ be the sequence generated by the EG-Anderson(1) algorithm, then the following statements hold.
\begin{enumerate}[{\rm (i)}]
\item (\textbf{Sequence convergence}) The sequence $\{x_k\}$ converges to a solution of ${\rm VI}(\Omega, H)$;

\

\item (\textbf{Best-iterate convergence rate}) $\min_{0\leq k\leq N}\|F_t(x_k)\|^2=O\left(\frac{1}{N}\right).$
\end{enumerate}
It can be seen that the EG-Anderson(1) algorithm can guarantee the sequence convergence as well as the EG method under the condition $tL<1$, and we will show that it is faster than the EG method by numerical experiments.
\end{remark}
\section{Numerical experiments}\label{s5}
In this section, we perform some numerical examples to compare the EG-Anderson(1) with Anderson(1) \cite{addition} and the EG method \cite{ref20}. 
All the codes were written in Matlab (R2023b) and run on a MacBook Air (16.00GB of RAM).

In the following numerical experiments, the stopping rule for Examples \ref{exam1}-\ref{exam4} is set by $$\|r(x_k)\| \coloneqq \|F_1(x_k)\| = \|x_k-P_{\Omega}(x_k-H(x_k))\|<10^{-8}$$ or the maximum iteration exceeds $10^4$ times. The parameters in the EG-Anderson(1) are set as follows
$$\omega=30, M=5000, \tau=0.6, \rho=0.8, \mu=0.5.$$
In the figures and tables of this section, 'Sec.' represents the CPU time in seconds and 'Iter.' represents the number of iterations.
Specifically, the numbers in parentheses represent the number of Anderson step \eqref{iteration3km} executed. Moreover, '$\backslash$' indicates that the number of iterations exceeds $10^4$, and the corresponding CPU time is not counted, represented by  $\textendash$.
Furthermore, the best performing algorithm in terms of the average number of iterations and CPU time is highlighted in bold for each combination of dimension $n$ and parameter $\gamma$.
\begin{example}\label{exam1}\cite{ref36}
Consider the Harker-Pang problem with linear mapping $H(x) \coloneqq Wx+w_0$, where $w_0\in\mathbb{R}^n$ and $$W \coloneqq A^{\rm T}A+S+D.$$
Here, $A$ is an $n\times n$ matrix, $S$ is an $n\times n$ skew-symmetric matrix and $D$ is an $n\times n$ diagonal matrix with nonnegative diagonal entries. Therefore, it follows that $W$ is positive semidefinite. Let the feasible set be $\Omega \coloneqq \{x\in\mathbb{R}^n:\textbf{0}\leq x\leq 20l\},$ where $l=(1,1,\ldots,1)^{{\rm T}}\in \mathbb{R}^n$.
It is clear that $H$ is monotone and Lipschitz continuous.
\end{example}
We can easily obtain that the Lipschitz constant of $H$ is $L = \|W\|$.
Applying the EG-Anderson(1) to this example, instead of using line search, we can do experiment with a constant stepsize $t$ that satisfies $tL<1$.

In the following experiments, we let $t=\frac{0.7}{L}$, and every entry of the skew-symmetric matrix $S$ is uniformly generated from $(-5, 5)$, and every diagonal entry of $D$ is uniformly generated from $(0, 2)$, and $A, w_0$ are randomly generated.

\begin{figure}[h]
\centerline{
 \subfigure[$n=500$]{\includegraphics[width=0.53 \textwidth]{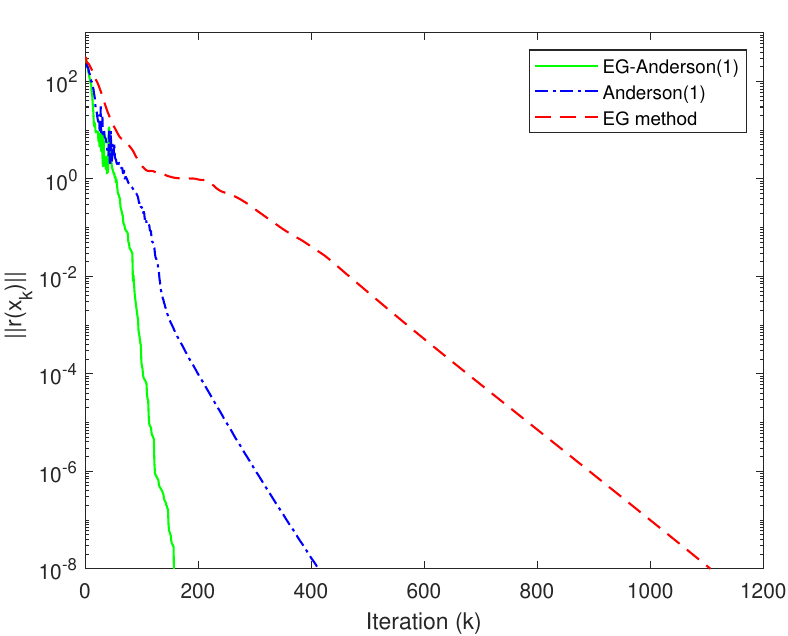}}
    \subfigure[$n=2000$]{\includegraphics[width=0.53 \textwidth]{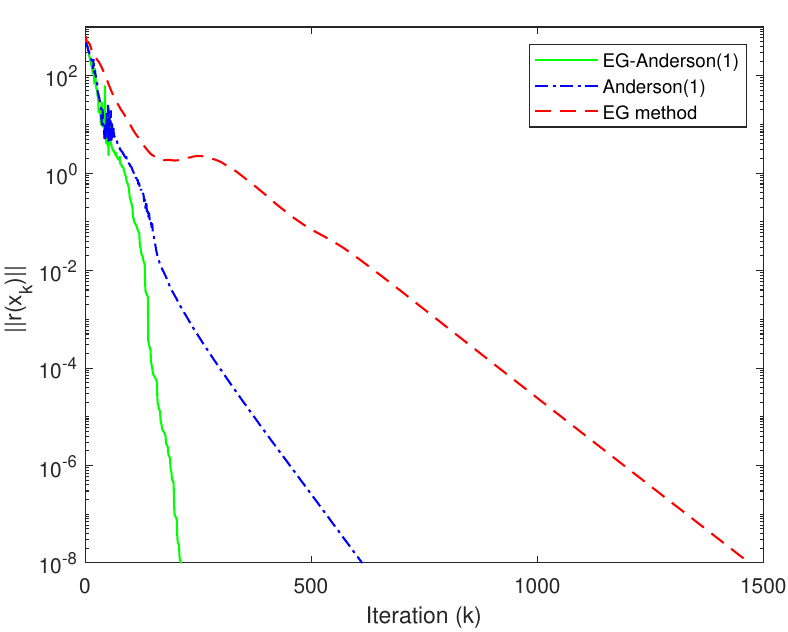}}
}
  \caption{Comparisons of the convergence behaviours of the EG-Anderson(1), Anderson(1) and the EG algorithms for Example \ref{exam1}}\label{fig1}
\end{figure}

Figure \ref{fig1} compares the decreasing on the residual function by the EG-Anderson(1), Anderson(1) and the EG algorithms at the same random initial point for Example \ref{exam1} with $n=500$ and $n=2000$, respectively. For different dimension $n$, Table \ref{ta1} illustrates the average number of iterations and CPU time of the corresponding experiments at ten random initial points, where we see that the superiorities of the EG-Anderson(1) over Anderson(1) and the EG algorithms gradually emerges as the dimension increases.

\begin{table}[h]
	\centering
	\label{table2}
	\begin{tabular}{c|ccccc}
		\toprule
		Sec.(avr)\\ Iter.(avr) &  & EG-Anderson(1) & Anderson(1) & EG \\
		\midrule
		\multirow{2}{1.5cm}{$n=100$} &  & 0.0038  & \textbf{0.0028}   & 0.0040   \\
		&  & \textbf{97.8} (96.6)   & 226.4  & 501.7   \\
		\midrule
		\multirow{2}{1.5cm}{$n=500$} &  & \textbf{0.0233}  & 0.0243   & 0.1171  \\
		& & \textbf{169} (166.4)  & 402.6   & 1103.8 \\
		\midrule
		\multirow{2}{1.5cm}{$n=1000$} &  & \textbf{0.0972}  & 0.1172   & 0.5764  \\
		& & \textbf{197.5} (195.6)   & 545.8   & 1417.8 \\
		\midrule
		\multirow{2}{1.5cm}{$n=2000$} & & \textbf{0.3462} & 0.4789   & 2.2768 \\
		& & \textbf{214.5} (211.6)   & 595.4   & 1465.2\\
		\midrule
		\multirow{2}{1.5cm}{$n=5000$} &  &  \textbf{2.2158}  & 2.7644   & 14.1645\\
		& & \textbf{244} (240.6)   & 600.9   & 1538 \\
		\midrule
		\multirow{2}{1.5cm}{$n=10000$} &  &  \textbf{9.6154}  & 10.5268   & 63.4244 \\
		& & \textbf{260.3} (256.7)   & 563.4   & 1685 \\
		\bottomrule
	\end{tabular}
	\caption{Comparisons of the three algorithms for Example \ref{exam1}}\label{ta1}
\end{table}

\begin{figure}
\centerline{
 \subfigure[$n=500$]{\includegraphics[width=0.53 \textwidth]{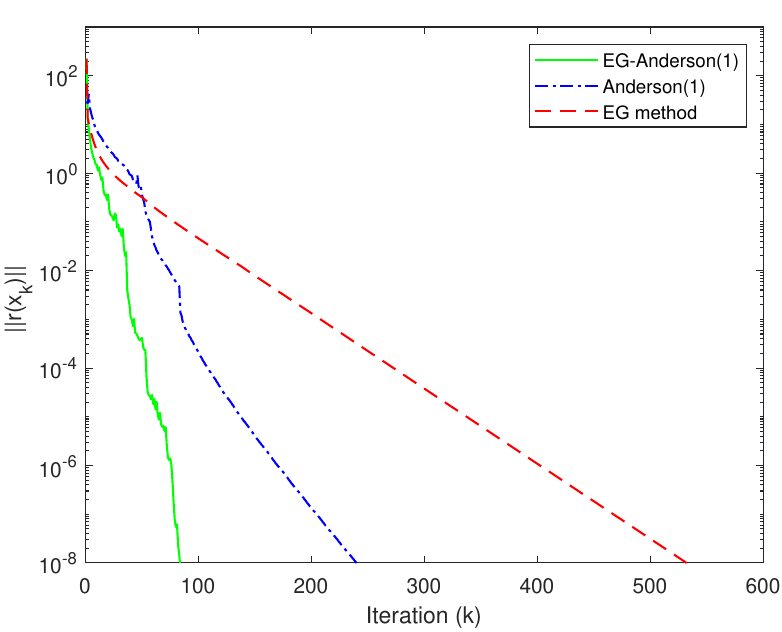}}
    \subfigure[$n=2000$]{\includegraphics[width=0.53 \textwidth]{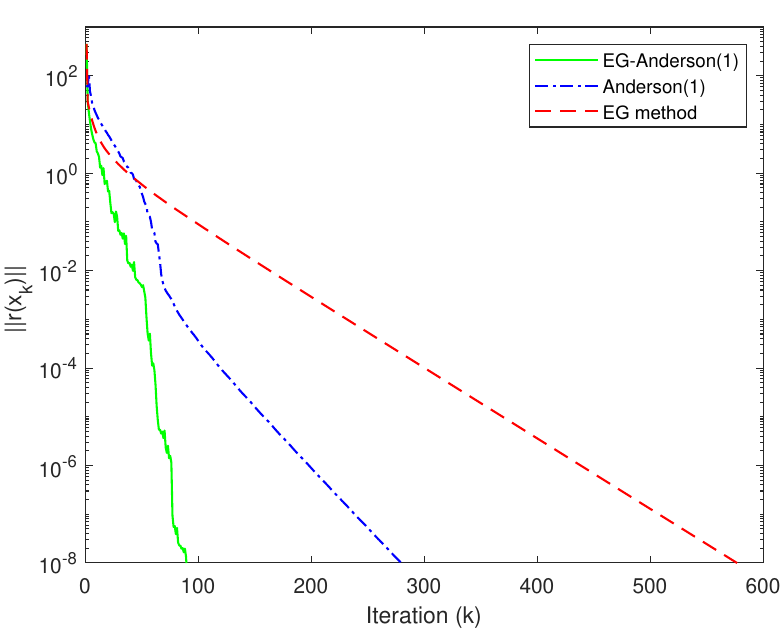}}
}
  \caption{Comparisons of the convergence behaviours of the EG-Anderson(1), Anderson(1) and the EG algorithms for Example \ref{exam2}}\label{fig21}
\end{figure}

\begin{example}\label{exam2}\cite{ref7}
Consider the quadratic fractional programming problem
\begin{eqnarray*}\label{sk}
\begin{array}{ll}
\mathrm{min}& \varphi(x) \coloneqq \frac{x^TQx+a^Tx+a_0}{b^Tx+b_0} \\
\mathrm{s.t.}& x\in \Omega \coloneqq \{x\in\mathbb{R}^n: 2l\leq x \leq 10l\}
\end{array}
\end{eqnarray*}
with $$Q\coloneqq Q_0^{{\rm T}}Q_0+I, a \coloneqq l+c, b \coloneqq l+d, a_0 \coloneqq 1+c_0, b_0 \coloneqq 1+d_0,$$
where $I\in\mathbb{R}^{n\times n}$ is the identity matrix, $l$ represents the vector that was defined in Example \ref{exam1} and $Q_0\in\mathbb{R}^{n\times n}$, $c,d \in \mathbb{R}^{n}$, $c_0, d_0\in\mathbb{R}$ are randomly generated from a uniform distribution.

It is easily verified that $\Omega\subseteq \{x\in\mathbb{R}^n:b^Tx+b_0>0\}$ and $Q$ is positive definite, and consequently $\varphi$ is pseudoconvex on $\Omega$.
Thus, $H(x)\coloneqq \nabla \varphi(x) $ in ${\rm VI}(\Omega, H)$ can be written in the following explicit form:$$H(x)=\frac{(b^Tx+b_0)(2Qx+a)-b(x^TQx+a^Tx+a_0)}{(b^Tx+b_0)^2}.$$
\end{example}
Let $\gamma=0.6$. Starting from the same random initial point, Figure \ref{fig21} shows the comparisons of the results obtained by the EG-Anderson(1), Anderson(1) and the EG algorithms for Example \ref{exam2} with $n=500$ and $n=2000$. Table \ref{ta2} presents the average number of iterations and CPU time for the experiments with more cases on $n$, conducted with ten different random initial points. Thus, the results indicate that the EG-Anderson(1) outperforms both Anderson(1) and the EG algorithms across various dimensions in terms of iterations and CPU time.

\begin{table}[h]
	\centering
	\label{table1}
	\begin{tabular}{c|ccccc}
		\toprule
		Sec.(avr)\\ Iter.(avr) &  & EG-Anderson(1) & Anderson(1) & EG \\
		\midrule
		\multirow{2}{1.5cm}{$n=100$} &  & 0.0051   & \textbf{0.0044}  & 0.0105  \\
		&   & \textbf{99} (99)  & 299.7   & 609  \\
		\midrule
		\multirow{2}{1.5cm}{$n=500$} &   & \textbf{0.0439}  & 0.0446  & 0.1819  \\
		&  & \textbf{93.2} (91.4)   & 238.5   & 531\\
		\midrule
		\multirow{2}{1.5cm}{$n=1000$} &  & \textbf{0.1631}  & 0.1937  & 0.7851  \\
		&  & \textbf{102.6} (101.3)   & 285.2   & 588.8 \\
		\midrule
		\multirow{2}{1.5cm}{$n=2000$} &   & \textbf{0.6839}  & 0.8337   & 3.4400  \\
		& & \textbf{101.2} (99.9)  & 285.3   & 576 \\
		\midrule
		\multirow{2}{1.5cm}{$n=5000$} &  & \textbf{5.5632}  & 6.9468   & 27.9809 \\
		& & \textbf{95.3} (93.7)  & 276   & 538 \\
		\midrule
		\multirow{2}{1.5cm}{$n=10000$} & & \textbf{18.7886}  & 20.7397   & 89.8027 \\
		& & \textbf{102.3} (100.9)   & 262   & 556 \\
		\bottomrule
	\end{tabular}
	\caption{Comparisons of the three algorithms for Example \ref{exam2}}\label{ta2}
\end{table}

\begin{example}\label{exam3}\cite{ref21}
Consider the following nonlinear complementarity problem (NCP)
\begin{equation}\label{NCP2}
H(x)\geq 0,~~~ x\geq 0, ~~~x^{\rm T}H(x)=0,
\end{equation}
where $H:\mathbb{R}^n\rightarrow\mathbb{R}^n$ is defined by
$$H(x)=(e^{-x^TUx}+\kappa)(Px+\iota).$$
Here $U$ is an $n\times n$ positive definite matrix, $P$ is an $n\times n$ positive semidefinite matrix, $\iota\in\mathbb{R}^n$ and $\kappa>0.$ It is also shown that $H$ is pseudomonotone \cite{ref21}.
Furthermore, the NCP in \eqref{NCP2} can be equivalently formulated as ${\rm VI}(\Omega, H)$ with $\Omega \coloneqq \{x\in\mathbb{R}^n: x\geq \textbf{0}\}$.
\end{example}

In numerical tests, we take $\kappa = 0.01, P=P_0^TP_0, U = U_0^TU_0$, where matrices $P_0,U_0\in\mathbb{R}^{n\times n}$ and vector $\iota\in\mathbb{R}^n$ are randomly generated as follows with an input integer $n$:
$$P_0={\rm randn}(n,n); U_0={\rm randn}(n,n);\hat{x}=\max\left(\textbf{0},{\rm randn}(n,1)\right);$$
$$ \iota = (-P*\hat{x}).*(\hat{x}>0)+(-P*\hat{x}+{\rm rand}(n,1)).*(\hat{x}==0).$$
From the above generation, we know that $\hat{x}$ is a solution of \eqref{NCP2}.

\begin{figure}
\centerline{
 \subfigure[$n=500$]{\includegraphics[width=0.53 \textwidth]{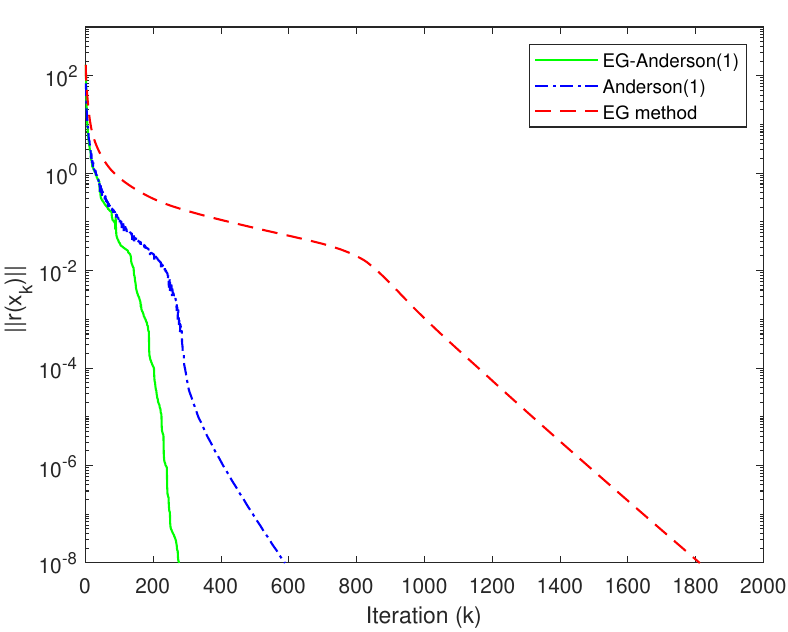}}
    \subfigure[$n=2000$]{\includegraphics[width=0.53 \textwidth]{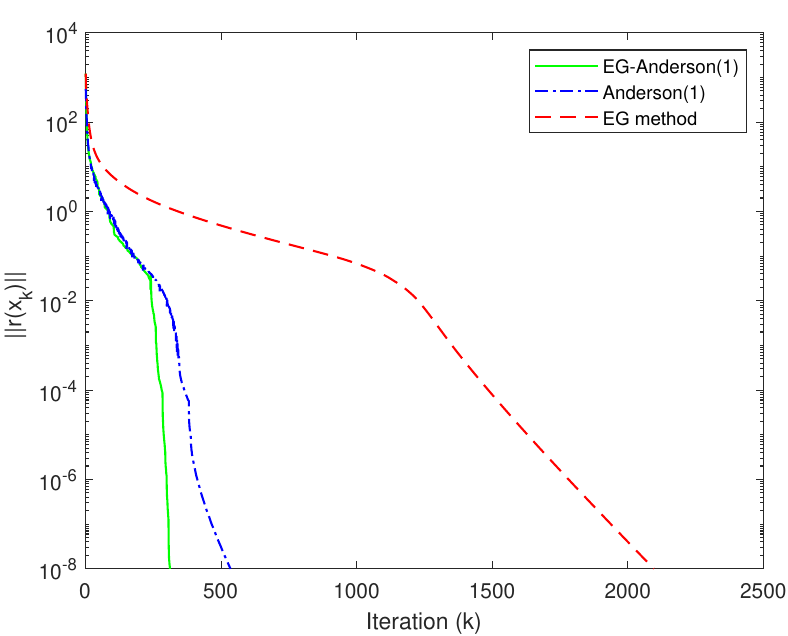}}
 }
  \caption{Comparisons of the convergence behaviours of the EG-Anderson(1), Anderson(1) and the EG algorithms for Example \ref{exam3}}\label{fig31}
\end{figure}

In Figure \ref{fig31}, we present a comparative analysis of the outcomes achieved by applying the EG-Anderson(1), Anderson(1) and the EG algorithms to Example \ref{exam3}. All three methods in Figure \ref{fig31} are initialized with the same random initial point. Table \ref{ta3} summarizes the average number of iterations and CPU time for each method, where the corresponding experiments are repeated ten times with different random initial points.

In all dimensions, the average number of iterations and CPU time of the EG-Anderson(1) are significantly less than those of Anderson(1) and the EG algorithms. The advantages of the EG-Anderson(1) become more pronounced as the dimension increases. For instance, in the case of $n=10000$, the EG-Anderson(1) achieves an average number of iterations of 552.2 and an average CPU time of 31.5962 seconds, significantly outperforming Anderson(1) and the EG algorithms.

\begin{table}[h]
	\centering
	\label{table1}
	\begin{tabular}{c|cccccc}
		\toprule
		Sec.(avr)\\ Iter.(avr) &  & EG-Anderson(1) & Anderson(1) & EG \\
		\midrule
		\multirow{2}{1.6cm}{$n=100$ $\gamma=0.27$} &  & 0.0041   & \textbf{0.0027}   & 0.0090  \\
		&   & \textbf{128.2} (125.4)  & 256.7   & 710.3  \\
		\midrule
		\multirow{2}{1.6cm}{$n=500$ $\gamma=0.03$} & &    \textbf{0.0405}  & 0.0424   & 0.2523  \\
		&  & \textbf{253.2} (249)   & 566.2  & 1754.4\\
		\midrule
		\multirow{2}{1.6cm}{$n=1000$ $\gamma=0.02$} & &   \textbf{0.1802} & 0.1823 & 1.0553  \\
		&   & \textbf{232.9} (228.8)   & 504   & 1513.1 \\
		\midrule
		\multirow{2}{1.6cm}{$n=2000$ $\gamma=0.008$} & &   \textbf{0.9466} & 0.9478  & 6.1649 \\
		&   & \textbf{321.2} (316.6)  & 622.3   & 2122.6 \\
		\midrule
		\multirow{2}{1.6cm}{$n=5000$ $\gamma=0.002$} &    & \textbf{7.3880}  & 8.3522  & 63.3128 \\
		&   & \textbf{485.2} (478.6)   & 1097   & 4146.3 \\
		\midrule
		\multirow{2}{1.6cm}{$n=10000$ $\gamma=0.0009$} &    & \textbf{31.5962}  & 33.8899   & 287.9682 \\
		&   & \textbf{552.2} (544.5)   & 1173.7   & 4836.6 \\
		\bottomrule
	\end{tabular}
	\caption{Comparisons of the three algorithms for Example \ref{exam3}}\label{ta3}
\end{table}

\begin{example}\label{exam4}
Consider the following partial differential equation (PDE) problem with free boundary
\begin{equation}\label{pde}
\begin{split}
 -\bigtriangleup u + \frac{9}{(1-p)^2}u^p+ \delta e^{-u}=c(\xi,\varsigma) \quad&{\rm in} ~\Lambda_+~~~
\\ u=0 ~~~~~~~~~~~~~~~~~\quad \quad \quad \quad \quad \quad &{\rm in} ~ \Lambda_0~~~
\\ u=|\nabla u|=0  ~~~~~~~\quad \quad \quad \quad &{\rm on}~~~~~ \Gamma~~~~
\\ u = v(\xi,\varsigma) ~~ ~~~~~~~~~~~~~\quad\quad\quad\quad &{\rm on} ~ \partial \Lambda,
\end{split}
\end{equation}
where $p\in(0,1)$, $\delta\geq 0$, $\Lambda=(0,1)\times(0,1)$, $\Lambda_+ = \{(\xi,\varsigma)\in\Lambda:u(\xi,\varsigma)>0\}$, $\Lambda_0 = \{(\xi,\varsigma)\in\Lambda: u(\xi,\varsigma)=0\}$ and $\Gamma=\partial \Lambda_0 = \partial \Lambda_+\cap\Lambda$ are unknown. Let $r^2=\xi^2+\varsigma^2$. We choose 
\begin{equation*}
	\begin{split} c(\xi,\varsigma)= D(r,p) R(r,p)+\delta e^{-R(r,p)}
	\end{split}
\end{equation*}
and $$v(\xi,\varsigma)=R(r,p),$$
where  $$D(r,p):=-\frac{3(3-p)\left[(3r-1)(1-p)+6r\right]}{r(1-p)^2(3r-1)^2}+\frac{27}{2(1-p)^2}\left(\frac{2}{3}\right)^p\left(\frac{3r-1}{2}\right)^{p-3}$$ and $R(r,p):=\left(\frac{3r-1}{2}\right)^{\frac{2}{1-p}}\max\left(0,r-\frac{1}{3}\right).$
\end{example}
Then problem \eqref{pde} has a solution as follows $$u(\xi,\varsigma)=R(r,p)=\left(\frac{3r-1}{2}\right)^{\frac{2}{1-p}}\max\left(0,r-\frac{1}{3}\right).$$
Dividing the interval (0,1) into $N$ subintervals of equal width $h$ provides mesh points $(\xi_i,\varsigma_j)$ where  $$\xi_i=ih, ~~i=0,1,\ldots,N$$$$\varsigma_j=jh,~~j=0,1,\ldots,N.$$
Using the five point finite difference method for the problem \eqref{pde} at grid $(\xi_i,\varsigma_j)$ gives
\begin{equation*}\label{dis}
-u_{i,j+1}-u_{i,j-1}+4u_{i,j}-u_{i+1,j}-u_{i-1,j}+\frac{9h^2}{(1-p)^2}u_{i,j}^p + h^2\delta e^{-u_{i,j}}=h^2c_{i,j}, ~~~~(\xi_i,\varsigma_j)\in\Lambda_+
\end{equation*}
and $$u_{i,j}=v_{i,j}, ~~~~ (\xi_i,\varsigma_j)\in \partial \Lambda.$$
Let $x:=(u_{1,1},u_{2,1},\ldots,u_{N-1,1},u_{1,2},\ldots,u_{N-1,N-1})^{\rm T}\in \mathbb{R}^{(N-1)^2}$ and $\tilde{c},\tilde{v}\in\mathbb{R}^{(N-1)^2}$ be the corresponding vectors transformed by $c_{i,j}$ and $v_{i,j}$.
Then, we obtain an NCP with $$H(x) \coloneqq Bx+E\max(0,x)^p+ Vf(x)+q,$$
where $B$ is a block tri-diagonal positive definite matrix of dimension $(N-1)^2\times (N-1)^2$,
$E$ and $V$ are both $(N-1)^2 \times (N-1)^2$ dimensional diagonal matrices with the diagonal elements being $\frac{9h^2}{(1-p)^2}$ and $\delta h^2$ respectively, $f(x):=\left(e^{-u_{1,1}}, e^{-u_{2,1}}, \ldots, e^{-u_{N-1,N-1}}\right)^{{\rm T}}\in\mathbb{R}^{(N-1)^2}$ and $q=-h^2\tilde{c}-\tilde{v}\in\mathbb{R}^{(N-1)^2}.$
Note that the dimension of the corresponding complementarity problem is $n=(N-1)^2$. Furthermore, the NCP can be equivalently formulated as ${\rm VI}(\Omega, H)$ with $\Omega:=\{x\in\mathbb{R}^n: x\geq\textbf{0}\}$.

\begin{figure}
\centering
\subfigure[$p=0.9$, $n=1600$, $\gamma=0.07$]{\includegraphics[width=0.48\textwidth]{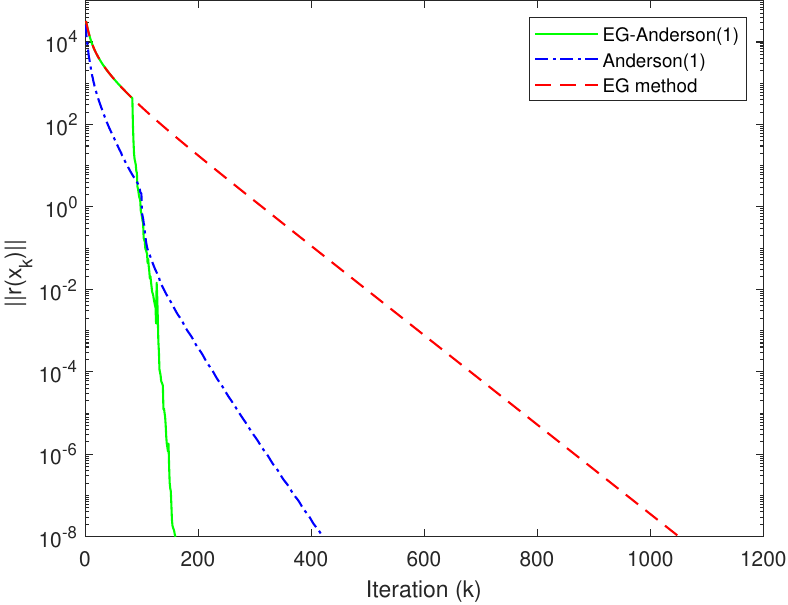}}
\hfill
\subfigure[$p=0.8$, $n=1600$, $\gamma=0.03$]{\includegraphics[width=0.48\textwidth]{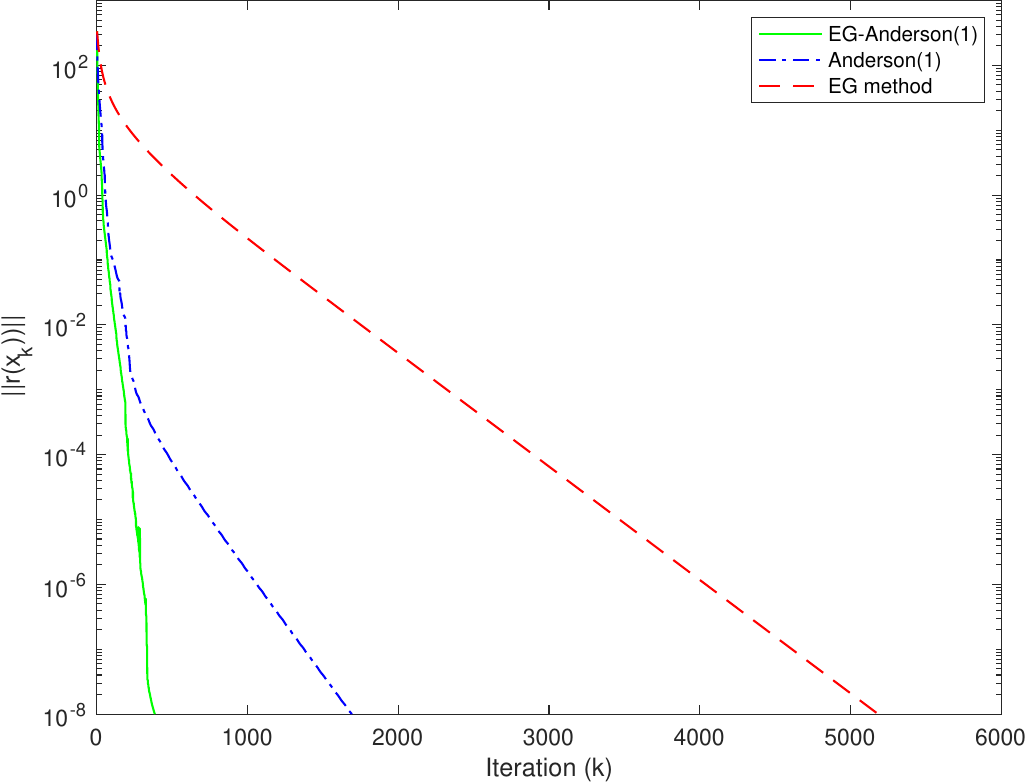}}

\vspace{1em}

\subfigure[$p=0.9, n=900, \gamma=0.39$]{\includegraphics[width=0.48\textwidth]{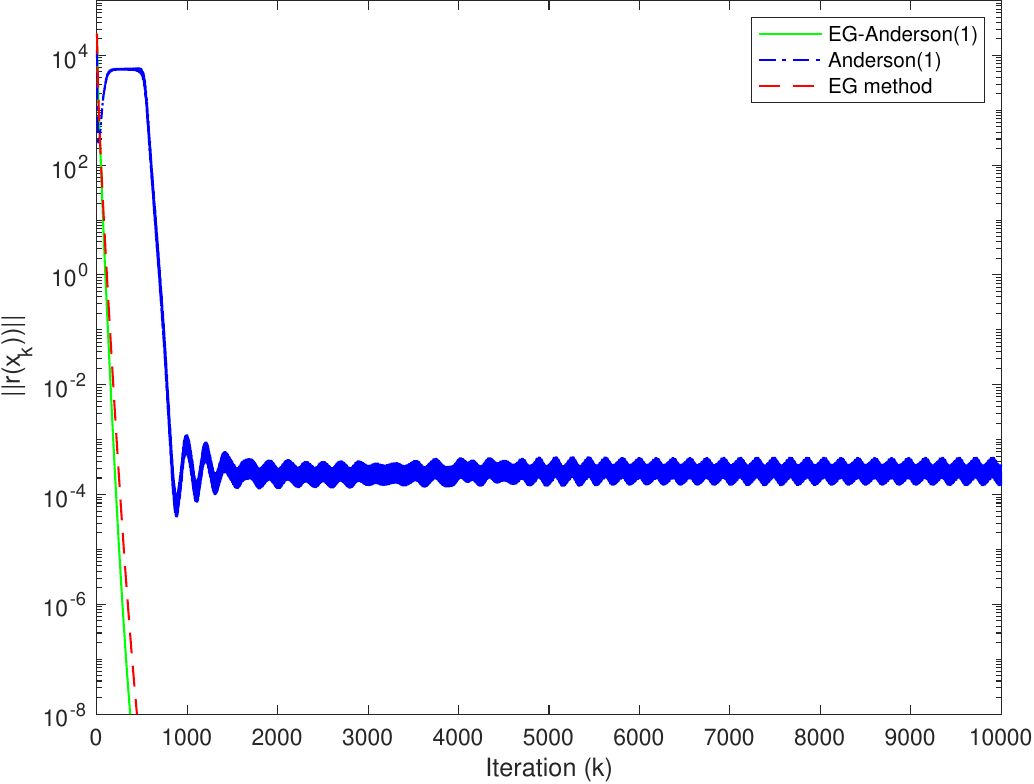}}
\hfill
\subfigure[$p=0.8, n=900, \gamma=0.39$]{\includegraphics[width=0.48\textwidth]{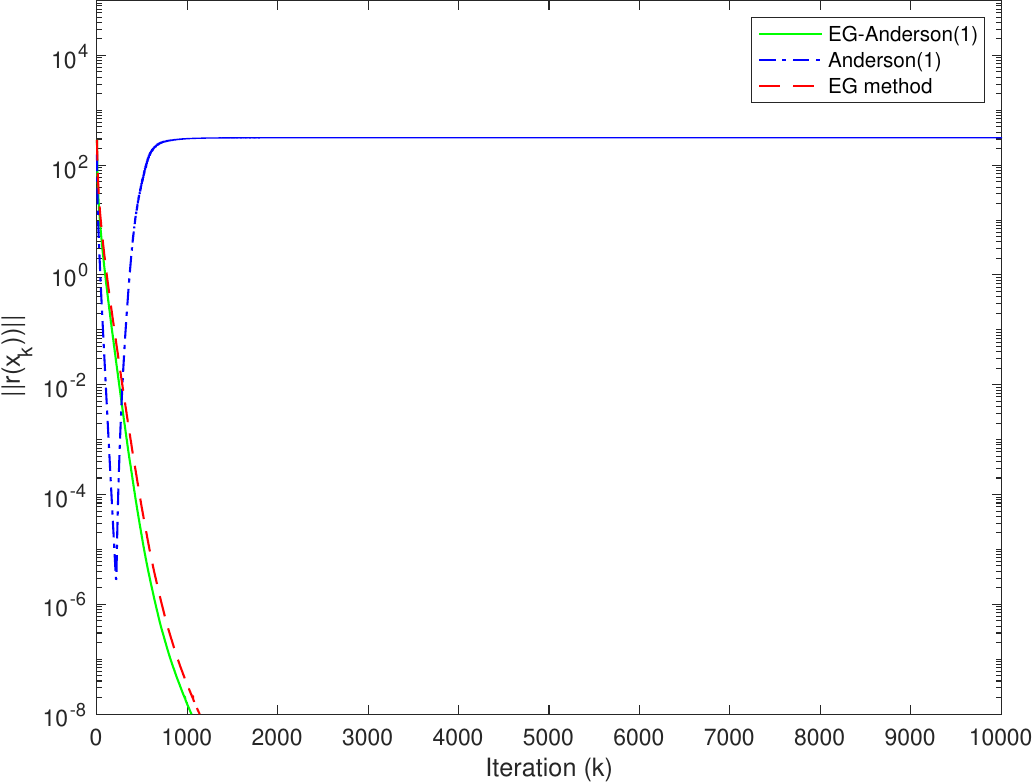}}

\caption{Comparisons of the convergence behaviours of the EG-Anderson(1), Anderson(1) and the EG algorithms for Example \ref{exam4}}
\label{fig411}
\end{figure}

In the experiments, let $\delta=1$. We set  $p=0.9$ and $p=0.8$, respectively.
The effectiveness of the EG-Anderson(1), Anderson(1) and the EG algorithms for Example \ref{exam4} with $n=900$ and $n=1600$ is compared in Figure \ref{fig411}. All three methods are started from the same random initial point for each case.

As can be seen from Table \ref{ta4}, the average number of iterations for the EG-Anderson(1) is consistently smaller than that of Anderson(1) and the EG algorithms in all dimensions. Additionally, the EG-Anderson(1) has a shorter CPU usage time, highlighting its higher computational efficiency.
In high-dimensional problems, such as $n=8100$ and $n=10000$, the EG-Anderson(1) outperforms the other two algorithms in terms of the average number of iterations and CPU time, demonstrating higher efficiency and performance.

Similar to Table \ref{ta4}, Table \ref{ta6} presents the average number of iterations and CPU time required for each algorithm in Example \ref{exam4}, but with a different value of $p$ (0.8 instead of 0.9).

\begin{table}[h]
	\centering
	\label{table1}
	\begin{tabular}{c|cccccc}
		\toprule
		Sec.(avr)\\ Iter.(avr) & & EG-Anderson(1) & Anderson(1) & EG \\
		\midrule
		\multirow{2}{1.5cm}{$n=100$ $\gamma=0.01$} &  & \textbf{0.0036}  & 0.0056  & 0.0135 \\
		&  & \textbf{64.7} (47.7)  & 280  & 530  \\
		\midrule
		\multirow{2}{1.5cm}{$n=900$ $\gamma=0.06$}& & 0.0165 & \textbf{0.0158}  & 0.0604 \\
		& & \textbf{121.1} (73.1)  & 300.3   & 656 \\
		\midrule
		\multirow{2}{1.5cm}{$n=1600$ $\gamma=0.07$} & & \textbf{0.0236}  & 0.0246 & 0.1022  \\
		& & \textbf{147.2} (79.2)   & 398.6  & 964 \\
		\midrule
		\multirow{2}{1.5cm}{$n=2500$ $\gamma=0.08$} & & \textbf{0.0349}  & 0.0412  & 0.1740 \\
		& & \textbf{186.9} (97.9)  & 526.7   & 1288 \\
		\midrule
		\multirow{2}{1.5cm}{$n=4900$ $\gamma=0.09$} &  & \textbf{0.0857}  & 0.1050  & 0.5587  \\
		& & \textbf{296.8} (156.8)  & 771.8   & 2178 \\
		\midrule
		\multirow{2}{1.5cm}{$n=8100$ $\gamma=0.1$} & & \textbf{0.1693}  & 0.1972   & 1.0951  \\
		& & \textbf{421.3} (228.3)  & 1064.9  & 3181 \\
		\midrule
		\multirow{2}{1.5cm}{$n=10000$ $\gamma=0.1$} & & \textbf{0.2659} & 0.3225   & 1.8438 \\
		& & \textbf{502.9} (276.9)   & 1270.7   & 3899 \\
		\bottomrule
	\end{tabular}
	\caption{Comparisons of the three algorithms for Example \ref{exam4} with $p=0.9$}\label{ta4}
\end{table}

\begin{table}[h]
	\centering
	\label{table1}
	\begin{tabular}{c|cccccc}
		\toprule
		Sec.(avr)\\ Iter.(avr) &  & EG-Anderson(1) & Anderson(1) & EG \\
		\midrule
		\multirow{2}{1.5cm}{$n=100$ $\gamma=0.005$} &  & 0.0115   & \textbf{0.0088}   & 0.0426 \\
		& & \textbf{371.2} (371.2)   & 774   & 2448.3  \\
		\midrule
		\multirow{2}{1.5cm}{$n=900$ $\gamma=0.02$} &  & \textbf{0.0383}  & 0.0683   & 0.3788  \\
		&  & \textbf{376.4} (376.4)   & 1478.4   & 4512.3\\
		\midrule
		\multirow{2}{1.5cm}{$n=1600$ $\gamma=0.03$} &  & \textbf{0.0843}   & 0.1161  & 0.8202 \\
		& & \textbf{375.1} (375.1) & 1427  & 5183.8 \\
		\midrule
		\multirow{2}{1.5cm}{$n=2500$ $\gamma=0.05$} &  & \textbf{0.0640}  & 0.1019  & 0.6231  \\
		&  & \textbf{397.9} (397.9)  & 1417.5   & 4761.1 \\
		\midrule
		\multirow{2}{1.5cm}{$n=4900$ $\gamma=0.06$} &  & \textbf{0.1355} & 0.2959  & 1.8412  \\
		&  & \textbf{474.4} (474.4)   & 2275.9  & 7554.5 \\
		\midrule
		\multirow{2}{1.5cm}{$n=8100$ $\gamma=0.07$} &  & \textbf{0.2472}  & 0.4831   & $\textendash$ \\ 
		&  & \textbf{639.4} (638.4)   & 2692.5   & $\backslash$ \\
		\midrule
		\multirow{2}{1.5cm}{$n=10000$ $\gamma=0.08$} &  & \textbf{0.3341}  & 0.6177   & $\textendash$ \\ 
		&  & \textbf{671.7} (669.7)   & 2660.8   & $\backslash$ \\
		\bottomrule
	\end{tabular}
	\caption{Comparisons of the three algorithms for Example \ref{exam4} with $p=0.8$}\label{ta6}
\end{table}

\begin{example}\label{ex5}
Consider the following linear complementarity problem (LCP)
\begin{equation}\label{LCP}
			\tilde{M}x+\tilde{q}\geq 0,~~~ x\geq 0, ~~~x^{\rm T}(\tilde{M}x+\tilde{q})=0,
\end{equation}
where $\tilde{M}$ is an $n\times n$ $P$-matrix and $\tilde{q}\in\mathbb{R}^n$ is a vector. Additionally, the LCP in \eqref{LCP} can be equivalently expressed as a nonmonotone ${\rm VI}(\Omega, H)$, where $H(x)=\tilde{M}x+\tilde{q}$ and $\Omega=\{x\in\mathbb{R}^n: x\geq \textbf{0}\}.$
\end{example}
Note that an $n\times n$ matrix $\tilde{M} = (m_{ij})$ is called a $P$-matrix if all principal minors of $\tilde{M}$ are positive.
We refer to Example 4.4 from \cite{pmatrix} to generate the matrix $\tilde{M}$ and vector $\tilde{q}$. First, we randomly generate a dense matrix $\tilde{A}\in\mathbb{R}^{n\times n}$ 
and a vector $\tilde{q}\in\mathbb{R}^n$, with elements uniformly distributed in the range $(-5,5)$. By applying the QR decomposition to $\tilde{A}$, we obtain an upper triangular matrix $\tilde{N}$. Next, we replace the diagonal elements of $\tilde{N}$ with their absolute values, resulting in a triangular matrix $\tilde{M}$ with positive diagonal entries. 
This ensures that $\tilde{M}$ is a $P$-matrix.

Below, we compare the performance of the EG-Anderson(1) algorithm, iPCA from \cite{new1}, IPA$_L$ from \cite{new3} and the semi-smooth Newton algorithm\footnote{The solver, developed by Y. Tassa, https://www.mathworks.com/matlabcentral/fileexchange/20952-lcp-mcp-solver-newton-based} from \cite{semi} for solving \eqref{LCP}. 
The stopping criterion for the experiments is the same as that for the semi-smooth Newton algorithm, specifically, $$0.5\left\|s_k-x_k-(\tilde{M}x_k+\tilde{q})\right\|^2<10^{-8}$$ or the maximum iteration exceeds $10^4$ times, where $s_k=((s_k)_1, (s_k)_2,\dots)^{\rm T}$ with
$(s_k)_i=\sqrt{(x_k)_i^2+(\tilde{M}x_k+\tilde{q})_i^2}$, $i=1,\dots, n$. Note that for an $n$-dimensional vector $z$, $(z)_i$ represents the $i$-th component of $z$.

We can find that the Lipschitz constant of $H$ is $\|\tilde{M}\|$. When we run the EG-Anderson(1) to this example, rather than employing a line search, we do experiment with a constant step size $t$ that meets the condition  $t\|\tilde{M}\|<1$.
In the following experiments, we let $t=\frac{0.7}{\|\tilde{M}\|}$.
In both the iPCA and the IPA$_L$, the step size is taken to be the same as that in the EG-Anderson(1) algorithm. In addition, we set the parameter $\gamma=1.5$ in the iPCA.

Table \ref{ta7} presents the average number of iterations and CPU time required by the four algorithms to meet the stopping criterion for \eqref{LCP} across different dimensions, starting from ten random initial points.
It can be observed that the EG-Anderson(1) algorithm consistently outperforms the other three algorithms in terms of CPU time across various dimensions.

\begin{table}[h]
	\centering
	\label{table7}
	\begin{tabular}{c|cccccc}
		\toprule
		Sec.(avr)\\ Iter.(avr) & & EG-Anderson(1) & semi-smooth Newton & iPCA & IPA$_L$ \\
		\midrule
		\multirow{2}{1.5cm}{$n=100$ } & & \textbf{0.0038} & 0.0237 & 0.0044 & 0.0215 \\
		& & 64 (63.2) & \textbf{29.4} & 149 & 281.5\\
		\midrule
		\multirow{2}{1.5cm}{$n=500$} & & \textbf{0.0188} & 0.3223 & 0.0464 & 0.4225 \\
		& & 125.2 (124.8) & \textbf{45.7} & 286.8 & 638.7\\
		\midrule
		\multirow{2}{1.5cm}{$n=1000$} & & \textbf{0.1372} & 1.3163 & 0.7139 & 19.1573 \\
		& & 293.6 ( 288.2) & \textbf{50.8} & 1221.2 & 2718.1 \\
		\midrule
		\multirow{2}{1.5cm}{$n=2000$} & & \textbf{0.5563} & 8.0195 & 2.3217 & 32.1018\\
		& & 306.4 (293.8) & \textbf{68.9} & 913.1 & 2048.6\\
		\midrule
		\multirow{2}{1.5cm}{$n=5000$ } & & \textbf{4.3059} & 17.1193 & 959.8955 & 854.2569 \\
		& & \textbf{408.3} (337.7) & 578.4 & 976.8 & 4994.4\\
		\bottomrule
	\end{tabular}
	\caption{Comparisons of the four algorithms for Example \ref{ex5}}\label{ta7}
\end{table}

\section{Conclusions}
This paper proposes an algorithm, called EG-Anderson(1) algorithm, for solving the pseudomonotone variational inequalities ${\rm VI}(\Omega, H)$. This algorithm is based on the EG method and Anderson acceleration.
Firstly, the global sequence convergence of the EG-Anderson(1) algorithm is proven without relying on the Lipschitz continuity and contractive condition that are required for the convergence analysis of the EG method and Anderson acceleration in prior research.
Moreover, when $H$ is locally Lipschitz continuous, the convergence rate of the residual function is analyzed and shown to be no worse than that of the EG method. Finally, the effectiveness of the EG-Anderson(1) algorithm has been validated through numerical experiments. The results demonstrate that it outperforms both Anderson(1) and the EG algorithms in terms of the number of iterations and CPU time, especially in the context of solving Harker-Pang problems, fractional programming problems, nonlinear complementarity problem, PDE problems with free boundary and linear complementarity problems.

\end{document}